\documentclass[11pt]{amsart}
\usepackage{amsmath,amssymb}
\usepackage{fullpage}
\usepackage[dvips]{graphicx}
\usepackage{psfrag}  
\def\COMMENT#1{}
\def\TASK#1{}
\def\noproof{{\unskip\nobreak\hfill\penalty50\hskip2em\hbox{}\nobreak\hfill%
        $\square$\parfillskip=0pt\finalhyphendemerits=0\par}\goodbreak}
\def\endproof{\noproof\bigskip}
\newdimen\margin   
\def\textno#1&#2\par{%
    \margin=\hsize
    \advance\margin by -4\parindent
           \setbox1=\hbox{\sl#1}%
    \ifdim\wd1 < \margin
       $$\box1\eqno#2$$%
    \else
       \bigbreak
       \hbox to \hsize{\indent$\vcenter{\advance\hsize by -3\parindent
       \sl\noindent#1}\hfil#2$}%
       \bigbreak
    \fi}
\def\proof{\removelastskip\penalty55\medskip\noindent{\bf Proof. }}

\def\eps{\varepsilon}

\newtheorem{firstthm}{Proposition}[section]
\newtheorem{thm}[firstthm]{Theorem}
\newtheorem{prop}[firstthm]{Proposition}
\newtheorem{lemma}[firstthm]{Lemma}

\newtheorem{conj}[firstthm]{Conjecture}
\newtheorem{claim}[firstthm]{Claim}
\newtheorem{fact}[firstthm]{Fact}
\newtheorem{ques}[firstthm]{Question}

\title{A degree sequence Hajnal--Szemer\'edi theorem}
\author{Andrew Treglown}

\begin{document}
\label{firstpage}
\date{\today}
\begin{abstract} 
We say that a graph $G$ has a perfect $H$-packing if there exists a set of vertex-disjoint 
copies of $H$ which cover all the vertices in $G$. The seminal Hajnal--Szemer\'edi theorem~\cite{hs} characterises the minimum degree
that ensures a graph $G$ contains a perfect $K_r$-packing.
Balogh, Kostochka and Treglown~\cite{bkt} proposed a degree sequence version of the Hajnal--Szemer\'edi theorem which, if true, gives a strengthening of the Hajnal--Szemer\'edi theorem. In this paper we prove this conjecture asymptotically. Another fundamental result in the area is the Alon--Yuster theorem~\cite{ay} which gives a minimum degree condition that ensures a graph contains a perfect $H$-packing for an \emph{arbitrary} graph $H$. We give a wide-reaching generalisation of this result by answering another conjecture of Balogh, Kostochka and Treglown~\cite{bkt} on the degree sequence of a graph that forces a perfect $H$-packing.
We also prove a degree sequence result concerning perfect transitive tournament packings in  directed graphs. The proofs blend together the regularity and absorbing methods. 
\end{abstract} 
\maketitle

\section{Introduction} 
\subsection{Perfect clique packings}
A \emph{perfect matching} in a graph $G$ is a collection of vertex-disjoint edges that together cover all the vertices of $G$.
A theorem of Tutte~\cite{tutte}  characterises those graphs that contain a perfect matching.
Moreover, the decision problem of whether a graph contains a perfect matching is polynomial time solvable~\cite{ed}.

A natural generalisation of the notion of a perfect matching is  a so-called perfect packing: Given two graphs $H$  and $G$, an \emph{$H$-packing} in $G$ 
is a collection of vertex-disjoint copies of $H$ in $G$. An
$H$-packing is called \emph{perfect} if it covers all the vertices of $G$.
Perfect $H$-packings are also referred to as \emph{$H$-factors} or \emph{perfect $H$-tilings}. Hell and Kirkpatrick~\cite{hell} showed that the decision problem of
whether a graph $G$ has a perfect $H$-packing is NP-complete precisely when $H$ has a
component consisting of at least $3$ vertices. Therefore, for such graphs $H$, it is unlikely that there is a complete characterisation
of those graphs containing a perfect $H$-packing.

The following classical result of Hajnal and Szemer\'edi~\cite{hs} characterises the minimum degree that ensures a graph contains a perfect $K_r$-packing. 
\begin{thm}[Hajnal and Szemer\'edi~\cite{hs}]\label{hs}
Every graph $G$ whose order $n$
is divisible by $r$ and whose minimum degree satisfies $\delta (G) \geq (1-1/r)n$ contains a perfect $K_r$-packing. 
\end{thm}
It is easy to see that the minimum degree condition here cannot be lowered. Earlier, Corr\'adi and Hajnal~\cite{corradi} proved Theorem~\ref{hs} in the
case when $r=3$. More recently, Kierstead and Kostochka~\cite{short} gave a short proof of the Hajnal--Szemer\'edi theorem. 

Over the last twenty years the Hajnal--Szemer\'edi theorem has been generalised in a number of directions. For example,
Kierstead and Kostochka~\cite{kier} proved an \emph{Ore-type} analogue of
the Hajnal--Szemer\'edi theorem: If $G$ is a graph whose order $n$ is divisible by $r$, then $G$ contains a perfect $K_r$-packing
provided that $d(x)+d(y) \geq 2(1-1/r)n-1$ for all non-adjacent $x \not = y \in V(G)$. Together with Balogh and Kostochka~\cite{bkt}, the
author characterised the edge density required to ensure a graph of given minimum degree contains a perfect $K_r$-packing. In~\cite{bkt} the following conjecture was also proposed.

\begin{conj}[Balogh, Kostochka and Treglown~\cite{bkt}]\label{conj1}
Let $n,r \in \mathbb N$ such that $r $ divides $n$. Suppose that $G$ is a graph on $n$ vertices with degree sequence $d_1 \leq 
\dots \leq d_n$ such that:
\begin{itemize}
\item[($\alpha$)] $d_i \geq (r-2)n/r+i $ for all $i < n/r$;
\item[($\beta$)] $d_{n/r+1} \geq (r-1)n/r$.
\end{itemize}
Then $G$ contains a perfect $K_r$-packing.
\end{conj}
Note that Conjecture~\ref{conj1}, if true, is much stronger than the Hajnal--Szemer\'edi theorem since the
degree condition allows for $n/r$ vertices to have degree less than $(r-1)n/r$. Moreover, the degree sequence condition in Conjecture~\ref{conj1} is best-possible. Indeed, examples in Section~4 in~\cite{bkt} show that one cannot replace ($\alpha$) with
$d_i \geq (r-2)n/r+i -1$ for a \emph{single} value of $i <n/r$ and ($\beta$) cannot be replaced with $d_{n/r+1} \geq (r-1)n/r -1$.
Chv\'atal's theorem on Hamilton cycles~\cite{ch} implies Conjecture~\ref{conj1} in the case when $r=2$.
In~\cite{bkt} it was shown that 
Conjecture~\ref{conj1} is true under the additional assumption that no vertex $x \in V(G)$ of degree less than $(r-1)n/r$ lies in a copy of $K_{r+1}$.\COMMENT{conjecture doesn't imply Ore-result for $r \geq 3$. Take sets of C vertices, 2n/3 and n/3-C. C set sees all but 2n/3 set. n/3 set independent. Add a few edges. Then satisfies Ore condition but not d.s}

In this paper we prove the following asymptotic version of Conjecture~\ref{conj1}.
\begin{thm}\label{mainthm}
Let $\gamma >0$ and $r \in \mathbb N$. Then there exists an integer $n_0 =n_0(\gamma ,r)$ such that the following holds. Suppose that 
$G$ is a graph on $n\geq n_0$ vertices where $r$ divides $n$ and with degree sequence $d_1 \leq 
\dots \leq d_n$ such that:
\begin{itemize}
\item  $d_i \geq (r-2)n/r+i +\gamma n  $ for all $i < n/r$.
\end{itemize}
Then $G$ contains a perfect $K_r$-packing. 
\end{thm}
Keevash and Knox~\cite{kk} also  have announced a proof of Theorem~\ref{mainthm} in the case when $r=3$. Conjecture~\ref{conj1} considers a `P\'osa-type' degree sequence condition: P\'osa's theorem~\cite{posa}
states that a graph~$G$ on~$n\ge 3$
vertices has a Hamilton cycle if its degree sequence $d_1\leq  \dots \leq d_n$ satisfies $d_i \geq i+1$ 
for all $i<(n-1)/2$ and if additionally $d_{\lceil n/2\rceil} \geq \lceil n/2\rceil$ when~$n$ is odd. The aforementioned result of Chv\'atal~\cite{ch} generalises P\'osa's theorem by characterising those degree sequences which force Hamiltonicity. It would be interesting to establish an analogous result for perfect $K_r$-packings. 

After this paper was submitted, the author and Staden~\cite{katherine} have used Theorem~\ref{mainthm} to prove a degree sequence version of P\'osa's conjecture on the square of a Hamilton cycle.

\subsection{General graph packings}
The following result of Alon and Yuster~\cite{ay} initiated the topic of `generalising' the Hajnal--Szemer\'edi theorem to perfect $H$-packings for an \emph{arbitrary} graph $H$.
\begin{thm}[Alon and Yuster~\cite{ay}]\label{ay1}
Suppose that $\gamma >0$ and $H$ is a graph with $\chi (H)=r$. Then there exists an integer $n_0=n_0 (\gamma ,H)$
such that the following holds. If $G$ is a graph whose order $n \geq n_0$ is divisible by $|H|$ and
$$\delta (G) \geq (1-1/r+\gamma )n$$
then $G$ contains a perfect $H$-packing.
\end{thm}
 Koml\'{o}s, S\'ark\"{ozy} and  Szemer\'{e}di~\cite{kss} proved that the  term $\gamma n$ in Theorem~\ref{ay1} can be replaced
 with a constant depending only on $H$. For some graphs $H$ the minimum degree condition in Theorem~\ref{ay1} is best-possible up to the  term $\gamma n$. However, for other graphs $H$ the minimum degree condition here can be substantially reduced. Indeed, K\"uhn and Osthus~\cite{kuhn2} characterised, up to an additive constant, the minimum degree which ensures a graph $G$ 
contains a perfect $H$-packing for an arbitrary graph $H$. This characterisation involves the so-called critical chromatic number of $H$ (see~\cite{kuhn2} for more details). K\"uhn, Osthus and Treglown~\cite{kotore}  characterised, asymptotically, the Ore-type degree condition that guarantees a graph $G$ 
contains a perfect $H$-packing for an arbitrary graph $H$.

In this paper we prove the following result, thereby answering another conjecture of Balogh, Kostochka and Treglown~\cite{bkt}.
\begin{thm}\label{conj2} Suppose that $\gamma >0$ and $H$ is a graph with $\chi (H)=r$. Then there exists an integer $n_0=n_0 (\gamma ,H)$
such that the following holds. If $G$ is a graph whose order $n \geq n_0$ is divisible by $|H|$, and whose degree
sequence $d_1\leq \dots \leq d_n$ satisfies
\begin{itemize}
\item  $d_i \geq (r-2)n/r+i +\gamma n  $ for all $i < n/r$,
\end{itemize}
then $G$ contains a perfect $H$-packing. 
\end{thm}
Note that Theorem~\ref{conj2} is a strong generalisation of the Alon--Yuster theorem.  Further, Theorem~\ref{mainthm} is a special case of Theorem~\ref{conj2} (so we do not prove the former directly). Previously, the author and~Knox~\cite{knox} proved Theorem~\ref{conj2} in the case when
$r=2$. (In fact, they proved a much more general result concerning embedding spanning bipartite graphs of small bandwidth.)

At first sight one may ask whether the  term $\gamma n$ in Theorem~\ref{conj2} can be replaced by a constant dependent on $H$. However, in Section~\ref{extremal} we show that for many graphs $H$ this is 
not the case. (In fact, for many graphs $H$ one cannot replace $\gamma n$ with $o(\sqrt{n})$.)

The proof of Theorem~\ref{conj2} splits into two main tasks: We use the regularity method to find an `almost' perfect $H$-packing in $G$ and find a so-called `absorbing set' that can be used to cover the remaining
vertices with disjoint copies of $H$ (see Section~\ref{subabs1} for the precise definition of such a set).

\subsection{Perfect packings in directed graphs}
Recently, there has been a focus on generalising the Hajnal--Szemer\'edi theorem to the directed graph setting (see for example~\cite{blm, cdkm, ckm, directed}). 
Here we give a degree sequence condition that forces a digraph to contain a perfect packing of transitive tournaments of a given size.
 Before we can state this result we require some notation and definitions.

Given a digraph $G$ on $n$ vertices and $x \in V(G)$, let $d^+_G(x)$ and $d^-_G(x)$ denote the out- and indegree of $x$ in $G$. We call $d^* _G (x) := \max \{ d^+ _G (x) , d^- _G (x) \}$ the \emph{dominant degree} of $x$ in $G$. The \emph{dominant degree sequence of $G$} is the sequence $d^* _1 \leq \dots \leq d^* _n$ of the dominant degrees of the vertices of $G$ ordered from smallest to largest.
Let $T_r$ denote the transitive tournament on $r$ vertices. 
\begin{thm}\label{mainthm2} Let $r \in \mathbb N$ and $\gamma >0$. Then there exists an integer $n_0=n_0 (\gamma ,r)$
such that the following holds. Let $G$ be a digraph on $n\geq n_0$ vertices where $r$ divides $n$. If $G$ has dominant degree sequence $d^* _1 \leq d^* _2 \leq \dots \leq d^* _n$ where
\begin{align*} d^* _i \geq (r-2)n/r+i +\gamma n \ \text{ for all } \ i<n/r
\end{align*}
then $G$ contains a perfect $T _r$-packing.
\end{thm}
Notice that Theorem~\ref{mainthm2} implies Theorem~\ref{mainthm}. (In particular, since Theorem~\ref{mainthm} is best-possible up to the  term $\gamma n$, so is Theorem~\ref{mainthm2}.)
 Further,
Theorem~\ref{mainthm2} complements a number of existing results on perfect $T_r$-packings: Let $n$ be divisible by $r$. Czygrinow, DeBiasio,  Kierstead and  Molla~\cite{cdkm} showed that a digraph on $n$
vertices contains a perfect $T_r$-packing provided that $G$ has minimum outdegree at least $(r-1)n/r$ or minimum total degree at least $2(r-1)n/r-1$. Amongst other results, in~\cite{directed} it was shown that if $G$ has minimum dominant
degree at least $(1-1/r+o(1))n$ then $G$ contains a perfect $T_r$-packing. (So Theorem~\ref{mainthm2} generalises this result.)
It would be interesting to establish whether the  term $\gamma n$ in Theorem~\ref{mainthm2} can be replaced by an additive constant. 

\subsection{Organisation of the paper} 
The paper is organised as follows. In the next section we introduce some notation and definitions. In Section~\ref{extremal} we give an
extremal example which shows that, in general, one cannot replace the  term $\gamma n$ in Theorem~\ref{conj2} with $o(\sqrt{n})$.
In Section~\ref{der} we state the main auxiliary results of the paper and derive Theorems~\ref{conj2} and~\ref{mainthm2} from them.
Szemer\'edi's Regularity lemma and related tools are introduced in Section~\ref{sec3}. We then prove an `almost' perfect packing result for Theorem~\ref{mainthm2} in Section~\ref{alsec}. In Section~\ref{overabs} we give an overview of our absorbing method and then apply it in Section~\ref{sec8} to prove  absorbing results for Theorems~\ref{conj2} and~\ref{mainthm2}.

\section{Notation and preliminaries}
\subsection{Definitions and notation}
Let $G$ be a (di)graph.
We write $V(G)$ for the vertex set of $G$, $E(G)$ for the edge set of $G$ and define
$|G|:=|V(G)|$ and $e(G):= |E(G)|$. Given a subset $X \subseteq V(G)$, we write $G[X]$ for the sub(di)graph of $G$ induced by $X$. 
We write $G\setminus X$ for the sub(di)graph of $G$ induced by $V(G)\setminus X$. Given a set $X \subseteq V(G)$ and a (di)graph $H$ on $|X|$ vertices we say that $X$ \emph{spans a copy of $H$ 
in $G$} if $G[X]$ contains a copy of $H$. In particular, this does not necessarily mean that
$X$ induces a copy of $H$ in $G$.

Suppose that $G$ is a graph. The degree of a vertex $x \in V(G)$ is denoted by $d_G(x)$ and its neighbourhood
by $N_G(x)$. Given a vertex $x \in V(G)$ and a set $Y \subseteq V(G)$ we write $d _G (x,Y)$ to denote the number of edges $xy$ where $y \in Y$.  Given disjoint vertex classes $A,B \subseteq V(G)$, $e_G (A,B)$ denotes the number of edges in $G$ with one endpoint in
$A$ and the other in $B$.

Given two vertices~$x$ and~$y$ of a digraph~$G$, we write~$xy$ for the edge directed from~$x$ to~$y$.
We denote by~$N^+ _G (x)$ and~$N^- _G (x)$ the out- and the inneighbourhood of~$x $
and by $d^+_G(x)$ and $d^-_G(x)$ its out- and indegree. 
We will write $N^+ (x)$ for example, if this is unambiguous. 
Given a vertex $x \in V(G)$ and a set $Y \subseteq V(G)$ we write $d^+ _G (x,Y)$ to denote the number of edges in $G$ with startpoint $x$ and
endpoint in $Y$. We define $d^- _G (x,Y)$ analogously.
Recall that  $d^* _G (x) := \max \{ d^+ _G (x) , d^- _G (x) \}$ is the \emph{dominant degree} of $x$ in $G$. If $d^* _G (x) = d^+ _G (x) = d^- _G (x) $ then we implicitly consider $d^* _G (x)$ to count the number
of edges sent out from $x$ in $G$. In particular, if $d^* _G (x) = d^+ _G (x) $ then we define $d^* _G (x,Y):=d^+ _G (x,Y)$. Otherwise, we set $d^* _G (x,Y):=d^- _G (x,Y)$.

Let $\mathcal H$ be a collection of (di)graphs and $G$ a (di)graph. An \emph{$\mathcal H$-packing} in $G$ is a collection of vertex-disjoint copies of elements
from $\mathcal H$ in $G$. Given an $\mathcal H$-packing $\mathcal M$, we write $V(\mathcal M)$ for the set of vertices covered by $\mathcal M$.

Recall that $T_r$ denotes the transitive tournament of $r$ vertices.
Given $1 \leq i \leq r$, we say a vertex $x \in V(T_r)$ is the \emph{$i$th vertex of $T_r$}
if $x$ has indegree $i-1$ and outdegree $r-i$ in $T_r$.

Given a (di)graph $G$ we let $G(t)$ denote the (di)graph obtained from $G$ by replacing each vertex $x \in V(G)$ with a set $V_x$ of $t$ vertices so that, for all $x,y \in V(G)$:
\begin{itemize}
\item If $xy \in E(G)$ then there are all possible edges in $G(t)$ with startpoint in $V_x$ and endpoint in $V_y$;
\item If $xy \not \in E(G)$ then there are no edges in $G(t)$ with startpoint in $V_x$ and endpoint in $V_y$.
\end{itemize}
We set $T_r ^t :=T_r (t)$ and $K_r ^t:=K _r (t)$.

Throughout the paper, we write $0<\alpha \ll \beta \ll \gamma$ to mean that we can choose the constants
$\alpha, \beta, \gamma$ from right to left. More
precisely, there are increasing functions $f$ and $g$ such that, given
$\gamma$, whenever we choose some $\beta \leq f(\gamma)$ and $\alpha \leq g(\beta)$, all
calculations needed in our proof are valid. 
Hierarchies of other lengths are defined in the obvious way.

\subsection{Absorbing sets}\label{subabs1}
Let $H$ be a (di)graph. Given a (di)graph $G$, a set $S \subseteq V(G)$ is called an \emph{$H$-absorbing set for $Q \subseteq V(G)$}, if both
$G[S]$ and $G[S\cup Q]$ contain perfect $H$-packings. In this case we say that \emph{$Q$ is $H$-absorbed by $S$}. Sometimes we will simply refer
to a set $S \subseteq V(G)$ as an $H$-absorbing set if \emph{there exists} a set $Q \subseteq V(G)$ that is $H$-absorbed by $S$.



\section{An extremal example for Theorem~\ref{conj2}}\label{extremal}
Let $K_{t_1, \dots , t_r}$ denote the complete $r$-partite graph with vertex classes of size $t_1, \dots , t_r$. The following result shows that for $H=K_{t_1, \dots , t_r}$ with $t_i \geq 2$ (for all $1 \leq i \leq r$) we cannot replace the  term $\gamma n$ in Theorem~\ref{conj2} with $\sqrt{n}/3r^2$. The construction given is a generalised version of an extremal graph from the arXiv version of~\cite{bkt}.

\begin{prop}\label{square} 
Let $r \geq 3$ and $H:=K_{t_1, \dots , t_r}$ with $t_i \geq 2$ (for all $1 \leq i \leq r$).
Let $n \in \mathbb N$ be sufficiently large so that $\sqrt{n}$ is an integer that is divisible by $6r^2|H|$. Set $C:=\sqrt{n}/3r^2$.
Then there exists a graph $G$ on $n$ vertices whose degree sequence $d_1 \leq \dots \leq d_n$ satisfies
$$d_i \geq (r-2)n/r +i + C \ \text{ for all } 1\leq i \leq \frac{n}{r} $$
but such that $G$ does not contain a perfect $H$-packing.
\end{prop}
\proof 
Let $G$ denote the graph on $n$ vertices consisting of $r$ vertex classes $V_1, \dots , V_r$ with $|V_1|=1$, $|V_2|= n/r+1+Cr$, $|V_3|= 2n/r -2 - 3 C$ and
$|V_i|=n/r- C$ if $4 \leq i \leq r$ and which contains the following edges:
\begin{itemize}
\item All possible edges with an endpoint in $V_3$ and the other endpoint in $V(G)\setminus V_1$. (In particular, $G[V_3]$ is complete.);
\item All edges with an endpoint in $V_2$ and the other endpoint in $V(G)\setminus V_2$;
\item All edges with an endpoint in $V_i$ and the other endpoint in $V(G)\setminus V_i$ for $4 \leq i \leq r$;
\item There are $\sqrt{n}/2$ vertex-disjoint stars in $V_2$, each of size $\lfloor 2|V_2|/\sqrt{n} \rfloor$,
$\lceil 2|V_2|/\sqrt{n} \rceil$, which cover all of $V_2$ (see Figure~1).
\end{itemize}
\begin{figure}\label{picture1}
\begin{center}\footnotesize
\includegraphics[width=0.5\columnwidth]{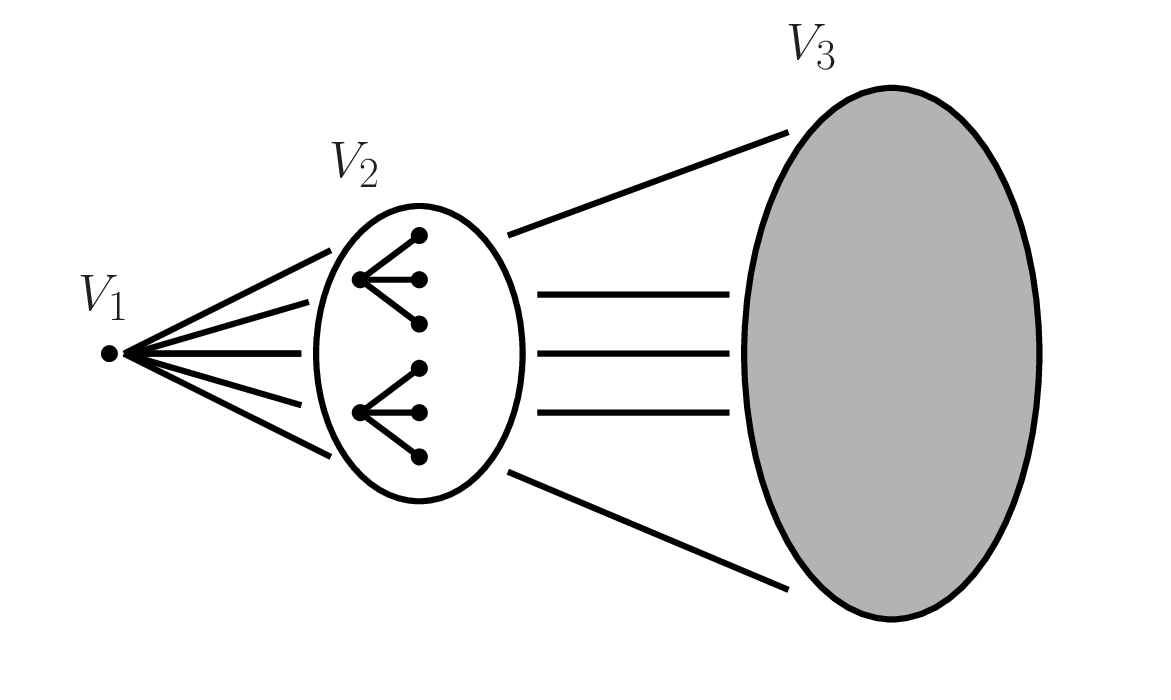}
\caption{An example of a graph $G$ from Proposition~\ref{square} in the case when $r=3$.}
\end{center}
\end{figure}
In particular, note that the vertex $v \in V_1$ sends all possible edges to $V(G) \setminus V_3$ but no edges to $V_3$.

Let $d_1 \leq \dots \leq d_n$ denote the degree sequence of $G$. 
Notice that every vertex in $V_i$ for $3 \leq i \leq r-3$ has degree at least $(r-1)n/r+C$. 
Note that $\lfloor 2|V_2|/\sqrt{n} \rfloor \geq 2 \sqrt{n}/r =6Cr$.
Thus, there are $\sqrt{n}/2$ vertices in $V_2$ of degree at least
$$(r-1)n/r-1-Cr +(6Cr-1)  \geq (r-1)n/r+C.$$
The remaining $n/r+1+Cr -\sqrt{n}/2 \leq n/r - \sqrt{n}/3 -1$ vertices in $V_2$ have degree at least
$$(r-1)n/r-Cr \geq (r-1)n/r  -\sqrt{n}/3 + C.$$
Since $d_G (v)\geq(r-2)n/r+1+C$ for the vertex $v \in V_1$ we have that
$d_i \geq (r-2)n/r +i + C$ for all  $1\leq i \leq {n}/{r} $.

Suppose that $v \in V_1$ lies in a copy $H'$ of $H$ in $G$.
Then by construction of $G$, two of the vertex classes $U_1,U_2$ of $H'$ must lie entirely in $V_2$. By definition of $H$, $H'[U_1\cup U_2]$ contains a path of length $3$. However, $G[V_2]$ does not contain a path of
length $3$, a contradiction. Thus, $v$ does not lie in a copy of $H$ and so $G$ does not contain a perfect $H$-packing.
\endproof
For many graphs $H$, K\"uhn and Osthus~\cite{kuhn2} showed that a minimum degree condition substantially smaller than that in the Alon--Yuster theorem ensures a graph $G$
contains a perfect $H$-packing.  For example, if $H=K_{t,t-1,\dots, t-1}$ for any $t \geq 3$ then they showed the so-called critical chromatic number of $H$, not $\chi (H)$, is the parameter that governs the minimum degree threshold that forces a perfect $H$-packing. Thus, it is interesting that for such graphs $H$ the degree sequence condition in Theorem~\ref{conj2}  is
`best-possible' up to the  term $\gamma n$. This raises the following natural question.
\begin{ques}
Are there graphs $H$ with $\chi (H)=r$ for which the degree sequence condition in Theorem~\ref{conj2} is `far' from tight?
\end{ques}

\section{Deriving Theorems~\ref{conj2} and~\ref{mainthm2}}\label{der}
\subsection{The auxiliary results for Theorem~\ref{mainthm2}}
The proof of Theorem~\ref{mainthm2} splits into two main tasks. Firstly, we construct a `small' $T_r$-absorbing set $M \subseteq V(G)$ with the property that both $G[M]$ and $G[M\cup Q]$ contain
perfect $T _r$-packings for \emph{any} `very small' set $Q \subseteq V(G)$ where $|Q|\in r \mathbb N$. Our second aim is to  find an `almost' perfect $T _r$-packing
in $G':=G\setminus M$. Together this ensures that $G$ contains a perfect $T _r$-packing. 
Indeed, suppose that $G'$ contains a $T _r$-packing $\mathcal M_1$ covering all but a very small set of vertices $Q$. Then by definition of $M$, $G[M\cup Q]$  contains a perfect $T _r$-packing $\mathcal M_2$. Thus, $\mathcal M_1 \cup \mathcal M_2$ is a perfect $T _r$-packing in $G$, as desired. 

The following result guarantees our desired absorbing set.
\begin{thm}\label{absthm}
Let $0 < 1/n \ll \nu \ll \eta  \ll 1/r$ where $n,r \in \mathbb N$ and $r \geq 2$. Suppose that $G$ is a digraph on $n$ vertices with dominant degree sequence $d^* _1 \leq d^* _2 \leq \dots \leq d^* _n$ where
\begin{align*} d^* _i \geq (r-2)n/r+i +\eta n \ \text{ for all } \ i<n/r.
\end{align*}
Then $V(G)$ contains a set $M$ so that $|M|\leq \nu n$ and $M$ is a $T_r$-absorbing set for any $W \subseteq V(G) \setminus M$ such that $|W| \in r \mathbb N$ and  $|W|\leq \nu ^3 n$.
\end{thm}
We prove Theorem~\ref{absthm} in Section~\ref{secab1}. In order to obtain our desired absorbing set $M$, we first show that $G$ contains many `connecting structures' of a certain type. Indeed, we show that there
is a constant $t \in \mathbb N$ such that given any $x,y \in V(G)$ there are `many' $t$-sets $X \subseteq V(G)$ so that both $G[X \cup \{x\}]$ and $G[X \cup \{y\}]$ contain perfect $T _r$-packings. A result of Lo and Markstr\"om~\cite{lo} then immediately implies the existence of the absorbing set $M$. A detailed exposition of our absorbing method is given in Section~\ref{overabs}.

The next result is used to obtain the almost perfect $T_r$-packing in $G'$.
\begin{thm}\label{approx} Let $n, r, h \in \mathbb N$ where $r \geq 2$ and $\eta >0$ such that $0<1/n \ll \eta \ll 1/r, 1/h$.  Let $G$ be a digraph on $n$ vertices with dominant degree sequence $d^* _1 \leq d^* _2 \leq \dots \leq d^* _n$ where
\begin{align*} d^* _i \geq (r-2)n/r+i +\eta n \ \text{ for all } \ i<n/r.
\end{align*}
Then $G$ contains a $T^h _r$-packing covering all but at most $\eta n$ vertices.
\end{thm}
Note that Theorem~\ref{approx} concerns $T^h_r$-packings rather than simply $T_r$-packings. (Recall that $T^h_r:=T_r (h)$.)
In Section~\ref{secapprox} we use the regularity method to prove Theorem~\ref{approx}. We give an outline of the proof  in Section~\ref{sketch}.

We now deduce Theorem~\ref{mainthm2} from Theorems~\ref{absthm} and~\ref{approx}.

\noindent
{\bf Proof of Theorem~\ref{mainthm2}.}
Define constants $\eta', \nu , \eta $ and  $n_0 \in \mathbb N$ such that
$0<1/n_0 \ll \eta ' \ll \nu \ll \eta \ll \gamma ,1/r.$
Let $G$ be a digraph on $n \geq n_0$ vertices as in the statement of the theorem. Apply Theorem~\ref{absthm} to $G$ with parameters $\nu, \eta$ to obtain a set $M$ so that $|M|\leq \nu n$ and $M$ is a $T_r$-absorbing set for any $W \subseteq V(G) \setminus M$ such that $|W| \in r \mathbb N$ and  $|W|\leq \nu ^3 n$. 

Set $G':=G\setminus M$ and $n':=|G'|$. It is easy to see that $G'$ has dominant degree sequence 
$d^* _{G',1} \leq d^* _{G',2} \leq \dots \leq d^* _{G',n'}$ where
\begin{align*} d^* _{G',i} \geq (r-2)n'/r+i +\eta ' n ' \ \text{ for all } \ i<n'/r.
\end{align*}
Apply Theorem~\ref{approx} to $G'$ with $\eta '$ playing the role of $\eta$. Then $G'$ contains a $T_r$-packing $\mathcal M_1$ covering all but a set $W$ of at most $\eta ' n' \leq \nu ^3 n$ vertices.
Since $n$ is divisible by $r$, $|W|$ is divisible by $r$. Hence, by definition of $M$, $G[M\cup W]$ contains a perfect $T_r$-packing $\mathcal M_2$. Then $\mathcal M_1 \cup \mathcal M_2$ is a perfect
$T_r$-packing in $G$, as desired.
\endproof

\subsection{The auxiliary results for  Theorem~\ref{conj2}.}\label{seccon}
The proof of Theorem~\ref{conj2} follows the same general strategy as the proof of Theorem~\ref{mainthm2}. The next result guarantees our desired absorbing set.
\begin{thm}\label{absorbG}
Let $n, r, h \in \mathbb N$ where $r \geq 2$ and $\nu, \eta >0$ such that $0<1/n \ll \nu \ll \eta \ll 1/r, 1/h$. 
Suppose that $H$ is a graph on $h$ vertices with $\chi (H)=r$.
 Let $G$ be a graph on $n$ vertices with degree sequence $d _1 \leq d _2 \leq \dots \leq d _n$ where
\begin{align*} d _i \geq (r-2)n/r+i +\eta n \ \text{ for all } \ i<n/r.
\end{align*}
Then $V(G)$ contains a set $M$ so that $|M|\leq \nu n$ and $M$ is an $H$-absorbing set for any $W \subseteq V(G) \setminus M$ such that $|W| \in h \mathbb N$ and  $|W|\leq \nu ^3 n$.
\end{thm}

Theorem~\ref{absorbG} is proved in Section~\ref{secabs2}. The next result is a simple consequence of Theorem~\ref{approx}.
\begin{thm}\label{approx2} Let $n, r, h \in \mathbb N$ where $r \geq 2$ and $\eta >0$ such that $0<1/n \ll \eta \ll 1/r, 1/h$. 
Suppose that $H$ is a graph on $h$ vertices with $\chi (H)=r$.
 Let $G$ be a graph on $n$ vertices with degree sequence $d _1 \leq d _2 \leq \dots \leq d _n$ where
\begin{align*} d _i \geq (r-2)n/r+i +\eta n \ \text{ for all } \ i<n/r.
\end{align*}
Then $G$ contains an $H $-packing covering all but at most $\eta n$ vertices.
\end{thm}
\proof
Let $G'$ denote the digraph obtained from $G$ by replacing each edge $xy \in E(G)$ with directed edges $xy,yx$. So $G'$ has dominant degree sequence 
$d^* _1 \leq d^* _2 \leq \dots \leq d^* _n$ where
\begin{align*} d^* _i \geq (r-2)n/r+i +\eta n \ \text{ for all } \ i<n/r.
\end{align*}
Thus, Theorem~\ref{approx} implies that $G'$ contains a $T^h _r$-packing $\mathcal M$ covering all but at most $\eta n$ vertices. By construction of $G'$, $\mathcal M$ corresponds to a $K^h _r$-packing in $G$.
Notice that $K^h _r$ contains a perfect $H$-packing. Therefore, $G$ contains an $H$-packing covering all but at most $\eta n$ vertices.
\endproof

Similarly to the proof of Theorem~\ref{mainthm2},
we now deduce Theorem~\ref{conj2} from Theorems~\ref{absorbG} and~\ref{approx2}.

\noindent
{\bf Proof of Theorem~\ref{conj2}.}
Define constants $\eta', \nu , \eta $ and  $n_0 \in \mathbb N$ such that
$0<1/n_0 \ll \eta ' \ll \nu \ll \eta \ll \gamma ,1/r,1/|H|.$ Set $h:=|H|$.
Let $G$ be a graph on $n \geq n_0$ vertices as in the statement of the theorem. Apply Theorem~\ref{absorbG} to $G$ with parameters $\nu, \eta$ to obtain a set $M$ so that $|M|\leq \nu n$ and $M$ is an $H$-absorbing set for any $W \subseteq V(G) \setminus M$ such that $|W| \in  h\mathbb N$ and  $|W|\leq \nu ^3 n$. 

Set $G':=G\setminus M$ and $n':=|G'|$. It is easy to see that $G'$ has  degree sequence 
$d _{G',1} \leq d _{G',2} \leq \dots \leq d _{G',n'}$ where
\begin{align*} d _{G',i} \geq (r-2)n'/r+i +\eta ' n ' \ \text{ for all } \ i<n'/r.
\end{align*}
Apply Theorem~\ref{approx2} to $G'$ with $\eta '$ playing the role of $\eta$. Then $G'$ contains an $H$-packing $\mathcal M_1$ covering all but a set $W$ of at most $\eta ' n' \leq \nu ^3 n$ vertices.
Since $n$ is divisible by $h$, $|W|$ is divisible by $h$. Hence, by definition of $M$, $G[M\cup W]$ contains a perfect $H$-packing $\mathcal M_2$. Then $\mathcal M_1 \cup \mathcal M_2$ is a perfect
$H$-packing in $G$, as desired.
\endproof

\section{The Regularity lemma and related tools}\label{sec3}
In the proofs of Theorems~\ref{approx} and~\ref{absorbG} we will apply
Szemer\'edi's Regularity lemma~\cite{reglem}. Before we  state it we need some more definitions. 
The \emph{density} of a bipartite graph $G=(A,B)$ with vertex classes~$A$ and~$B$ is defined to be 
$$d_G (A,B):=\frac{e_G(A,B)}{|A||B|}.$$ 
We will write $d(A,B)$ if this is unambiguous. Given any $\varepsilon >0$ we say that~$G$
is {\it $\varepsilon$-regular} if for all $X\subseteq A$ and $Y \subseteq B$ with $|X|>\varepsilon 
|A|$ and $|Y|> \varepsilon |B|$ we have that $|d(X,Y)-d(A,B)|<\varepsilon$.

Given disjoint vertex sets~$A$ and~$B$ in a graph~$G$, we write $(A,B)_G$
for the induced bipartite subgraph of~$G$ whose vertex classes are~$A$ and~$B$.
If $G$ is a digraph, we write $(A,B)_G$ for the oriented bipartite subgraph of~$G$
whose vertex classes are~$A$ and~$B$ and whose edges are all the edges from~$A$ to~$B$ in~$G$.
We say $(A,B)_G$ is {\it $\varepsilon$-regular and has density~$d$} if the
underlying bipartite graph of $(A,B)_G$ is $\varepsilon$-regular and has density~$d$.
(Note that the ordering of the pair $(A,B)$ is important.)

The Diregularity lemma is a variant of the Regularity lemma for digraphs due to Alon
and Shapira~\cite{alon}. Its proof is similar to the undirected version.
We will use the following degree form of the Regularity lemma which is stated in terms of both graphs and digraphs.  It is derived from the Diregularity lemma.
(See for example~\cite{ayphd} for a proof of the directed version and~\cite{survey} for a sketch of the undirected version.) 

\begin{lemma}[Degree form of the Regularity lemma]\label{dilemma}
For every $\eps\in (0,1)$ and every integer~$M'$ there are integers~$M$ and~$n_0$
such that if~$G$ is a graph or digraph on $n\ge n_0$ vertices and
$d\in[0,1]$ is any real number, then there is a partition of the vertices of~$G$ into
$V_0,V_1,\ldots,V_k$ and a spanning sub(di)graph~$G'$ of~$G$ such that the following holds:
\begin{itemize}
\item $M'\le k\leq M$;
\item $|V_0|\leq \varepsilon n$;
\item $|V_1|=\dots=|V_k|=:m$;
\item If $G$ is a graph then $d_{G'}(x)>d_G(x)-(d+\varepsilon)n$ for all vertices $x\in V(G)$;
\item If $G$ is a digraph then $d^+_{G'}(x)>d^+_G(x)-(d+\varepsilon)n$  and  $d^-_{G'}(x)>d^-_G(x)-(d+\varepsilon)n$ for all vertices $x\in V(G)$;
\item  $e(G'[V_i])=0$ for all $i \geq 1$;	
\item For all $1\leq i,j\leq k$ with $i\neq j$ the pair $(V_i,V_j)_{G'}$ is $\varepsilon$-regular
and has density either~$0$ or density at least~$d$.
\end{itemize}
\end{lemma}
We call $V_1, \dots, V_k$ \emph{clusters}, $V_0$ the \emph{exceptional set} and the
vertices in~$V_0$ \emph{exceptional vertices}. We refer to~$G'$ as the \emph{pure (di)graph}.
The \emph{reduced (di)graph~$R$ of~$G$ with parameters $\eps$, $d$ and~$M'$} is the (di)graph whose 
vertex set is $V_1, \dots , V_k$ and which $V_i V_j$ is an edge precisely when $(V_i,V_j)_{G'}$
is $\varepsilon$-regular and has density at least~$d$.

Roughly speaking, the next result states that the reduced digraph $R$ of a digraph $G$ as in Theorem~\ref{approx} `inherits' the dominant degree sequence of $G$. 
\begin{lemma}\label{inherit}
Let $M',n_0, r \in \mathbb N$ where $r \geq 2$ and let $\eps, d , \eta$ be positive constants such that
$1/n_0 \ll 1/M' \ll \eps \ll d \ll \eta <1$. Let $G$ be a digraph on $n \geq n_0$ vertices with dominant degree sequence $d^* _1 \leq d^* _2 \leq \dots \leq d^* _n$ 
such that
\begin{align}\label{ds1} d^* _{i} \geq (r-2)n/r+i+\eta n  \ \text{ for all } \ i<n/r.
\end{align}
Let $R$ be the reduced digraph of $G$ with parameters $\eps, d$ and $M'$ and set $k:=|R|$. Then $R$ has 
dominant degree sequence $d^* _{R,1} \leq d^* _{R,2} \leq \dots \leq d^* _{R,k}$ 
such that
\begin{align}\label{ds2} d^* _{R,i} \geq (r-2)k/r+i+\eta k/2  \ \text{ for all } \ i<k/r.
\end{align}
\end{lemma}
\proof
Let $G'$ denote the pure digraph of $G$ and let $V_1, \dots , V_k$ denote the clusters of $G$ and $V_0$ the exceptional set. Set $m:= |V_1|=\dots =|V_k|$. We may assume that
$d^* _R (V_1) \leq d^* _R (V_2) \leq \dots \leq d^* _R (V_k)$. Consider any $i <k/r$. Set $S:= \bigcup _{1 \leq j \leq i} V_j$. So $|S|=mi <mk/r\leq n/r$. Thus, by (\ref{ds1}) there exists a vertex
$x \in S$ such that $d^* _G (x) \geq d^* _{mi} \geq (r-2)n/r+mi+\eta n$.
Suppose that $x \in V_j$ where $1\leq j \leq i$. Since $km\leq n$, Lemma~\ref{dilemma} implies that
$$d^* _R (V_j) \geq (d^* _{G'} (x) -|V_0|)/m \geq ((r-2)n/r+mi +\eta n -(d+2\eps)n)/m \geq (r-2)k/r+i +\eta k/2.$$
Since $d^* _{R,i} =d^* _R (V_i) \geq d^* _R (V_j)$ this proves that (\ref{ds2}) is satisfied, as desired.
\endproof
An analogous proof yields the following version of Lemma~\ref{inherit} for graphs.
\begin{lemma}\label{inherit2}
Let $M',n_0, r \in \mathbb N$ where $r \geq 2$ and let $\eps, d , \eta$ be positive constants such that
$1/n_0 \ll 1/M' \ll \eps \ll d \ll \eta <1$. Let $G$ be a graph on $n \geq n_0$ vertices with degree sequence $d _1 \leq d _2 \leq \dots \leq d _n$ 
such that
\begin{align*} d _{i} \geq (r-2)n/r+i+\eta n  \ \text{ for all } \ i<n/r.
\end{align*}
Let $R$ be the reduced graph of $G$ with parameters $\eps, d$ and $M'$ and set $k:=|R|$. Then $R$ has 
degree sequence $d_{R,1} \leq d _{R,2} \leq \dots \leq d _{R,k}$ 
such that
\begin{align*} d _{R,i} \geq (r-2)k/r+i+\eta k/2  \ \text{ for all } \ i<k/r.
\end{align*}
\end{lemma}

The next well-known observation (see~\cite{ko} for example) states that a large subgraph of a  regular pair is also regular.
\begin{lemma}\label{slice}
Let $0< \eps < \alpha$ and $\eps ':= \max \{ \eps /\alpha , 2\eps \}$. Let $(A,B)$ be an $\eps$-regular pair of density $d$. Suppose $A' \subseteq A$ and $B' \subseteq B'$ where $|A'|\geq \alpha |A|$ and 
$|B'|\geq \alpha |B|$. Then $(A',B')$ is an $\eps '$-regular pair with density $d'$ where $|d'-d|<\eps$.
\end{lemma}

The following  result will be applied in the proof of Theorem~\ref{approx} to convert an almost perfect $T_r$-packing in a reduced digraph into an almost perfect $T^h _r$-packing in the original digraph $G$. It is (for example) a special case of Corollary 2.3 in~\cite{alony}.
\begin{lemma}\label{red}
Let $\eps, d>0$ and $m,r, h \in \mathbb N$ such that $0<1/m \ll \eps \ll d \ll 1/r, 1/h$. Let $H$ be a graph obtained from $K_r$ by replacing every vertex of $K_r$ with $m$ vertices and replacing each edge of $K_r$ with an $\eps ^2$-regular pair of density at least $d$. Then $H$
contains a $K^h _r$-packing covering all but at most $\eps m r$ vertices.
\end{lemma}

In the proof of Theorem~\ref{absorbG} we will apply the following well-known counting lemma (often called the Key lemma) from~\cite{simo}.
\begin{lemma}[Counting lemma~\cite{simo}]\label{simo}
Suppose that $0<\eps < d$, that $m,t \in \mathbb N$ and that $R$ is a graph with $V(R)=\{v_1, \dots , v_k\}$. Construct a graph $G$ by replacing every vertex $v_i \in V(R)$ by a set $V_i$ of $m$ vertices, and replacing
the edges of $R$ with $\eps$-regular pairs of density at least $d$. For each $v_i \in V(R)$ let $U_i$ denote the set of $t$ vertices in $R(t)$ corresponding to $v_i$. 
Let $H$ be a subgraph of $R(t)$ on $h$ vertices and maximum degree $\Delta \in \mathbb N$. Set $\delta :=d-\eps$ and $\eps _0 := \delta ^\Delta /(2+\Delta)$. If $\eps \leq \eps _0$ and $t-1 \leq \eps _0 m$ then
there are at least
$$ (\eps _0 m) ^h \text{ labelled copies of $H$ in $G$ }$$
so that if $x \in V(H)$ lies in $U_i$ then $x$ is embedded into $V_i$ in $G$.
\end{lemma}

\section{Almost perfect packings in digraphs}\label{alsec}
\subsection{Outline of the proof of Theorem~\ref{approx}}\label{sketch}
In Section~\ref{secapprox} we prove Theorem~\ref{approx}. The idea of the proof is as follows: We apply the Regularity lemma (Lemma~\ref{dilemma}) to obtain a reduced digraph $R$ of $G$. If $R$ contains an
almost perfect $T_r$-packing then by applying Lemma~\ref{red} we conclude that $G$ contains an almost perfect $T^h _r$-packing, as required. Otherwise, suppose that the largest $T_r$-packing  in $R$
covers precisely $d\leq (1-o(1))|R|$ vertices. We then show that there is a $\{T_r,T_{r+1}\}$-packing in $R$ covering substantially more than $d$ vertices (see Lemma~\ref{expand}). Since $T_r (r)$ and $T_{r+1} (r)$ both contain
perfect $T_r$-packings, this implies that the blow-up $R(r)$ of $R$ contains a $T_r$-packing covering substantially more than $dr$ vertices. Thus, crucially, the largest $T_r$-packing in $R(r)$ covers a higher proportion of vertices than the largest $T_r$-packing in $R$. By repeating this argument, we obtain a blow-up $R'$ of $R$ that contains an almost perfect $T_r$-packing. Using Lemma~\ref{red} we then
show that this implies that $G$ contains an almost perfect $T^h _r$-packing, as desired. Arguments of a similar flavour were applied in~\cite{ko, hlad}.

\subsection{Tools for the proof of Theorem~\ref{approx}}
In this section we  give a number of tools that will be used in the proof of Theorem~\ref{approx}.
The following trivial observation will be used in the proof  to convert a $\{T_r,T_{r+1}\}$-packing in a (reduced) digraph $R$ into a $T_r$-packing in the blow-up $R(r)$ of $R$.
\begin{fact}\label{fact2} Suppose that $r ,t\in \mathbb N$ such that $r$ divides $t$. Then both $T_r (t)=T^t_r$ and $T_{r+1} (t)= T^{t} _{r+1}$ contain  perfect $T_r$-packings.
\end{fact}
The next result will be used to show that a blow-up of a reduced digraph $R$ `inherits' the dominant degree sequence of $R$.
\begin{prop}\label{fact3}
Let $n,r,t \in \mathbb N$ and $\gamma >0$ where $n>r$. Suppose that $H$ is a digraph on $n$ vertices with dominant degree sequence $d^* _{H,1} \leq \dots \leq d^* _{H,n}$ such that
\begin{align}\label{lafact} d^* _{H,i} \geq (r-2)n/r+i +\gamma n \ \text{ for all } \ i<n/r.
\end{align}
Then 
$H':=H(t)$ has dominant degree sequence $d^* _{H',1} \leq \dots \leq d^* _{H',nt}$ such that
\begin{align*} d^* _{H',i} \geq (r-2)nt/r+i +(\gamma n-1)t \ \text{ for all } \ i<nt/r.
\end{align*}
\end{prop}
\proof
Given any $1 \leq j \leq nt$, by definition of $H'$ we have that
\begin{align}\label{compare}
d^* _{H',j} = t \times d^*_{H, \lceil j/t \rceil}.
\end{align}
Suppose that $j \leq nt/r -t$. Then $\lceil j/t \rceil< n/r$, so (\ref{lafact}) and (\ref{compare}) imply that
$$d^* _{H',j} \geq (r-2)nt/r +\lceil j/t \rceil  t +\gamma nt \geq (r-2)nt/r +j   +\gamma nt .$$
In particular, this implies that for any $i < nt/r$,
$$d^* _{H',i}  \geq (r-2)nt/r +(i-t) +\gamma nt = (r-2)nt/r+i +(\gamma n-1)t ,$$ as desired.
\endproof

\begin{prop}\label{fact}
Let $n,r  \in \mathbb N$ such that $n \geq r \geq 2$.  Let $G$ be a digraph on $n$ vertices with dominant degree sequence $d^* _1 \leq d^* _2 \leq \dots \leq d^* _n$ where
\begin{align}\label{min*} d^* _{\lceil {n/r} \rceil } \geq (1-1/r)n  .
\end{align}
Then $G$ contains a copy of $T_{r}$.
\end{prop}
\proof
Let $r' \in \mathbb N$. Suppose that $T'$ is a copy of $T_{r'}$ in $G$.
Let $V(T')=\{x_1, \dots , x_{r'}\}$ where $x_i$ plays the role of the $i$th vertex
of $T_{r'}$. We say that $T'$ is \emph{consistent} if there exists $0\leq s' \leq r'$
such that
\begin{itemize}
\item $d_G ^+ (x_i) \geq  ( 1- {1}/{r}  ) n$ for all $i \leq s'$;
\item $d_G ^- (x_i) \geq  ( 1- {1}/r ) n$ for all $i >s'$.
\end{itemize}
We call $s'$ a \emph{turning point of $T'$}. (Note that
$T'$ could have more than one turning point.)

(\ref{min*}) implies that there is a copy of $T_1$ in $G$ that is consistent.
Suppose that, for some $1 \leq r' < r$, we have found a consistent copy $T'$ of $T_{r'}$ in $G$.
As before, let $V(T')=\{x_1, \dots , x_{r'}\}$ where $x_i$ plays the role of the $i$th vertex
of $T_{r'}$ and let $s'$ denote a turning point of $T'$. 
Set $N':=\bigcap _{i\leq s'} N^+ _G(x_i) \cap \bigcap _{i>s'} N^- _G (x_i).$
Since $T'$ is consistent with turning point $s'$ and $r'\leq r-1$,
$$|N'|\geq \left (1 -\frac{r'}{r}  \right )n \geq \frac{n}{r} .$$
In particular, (\ref{min*}) implies there is a vertex $x \in N'$ such that
$d^* _G(x) \geq (1-1/r)n$. Then $V(T') \cup \{x\}$ spans a consistent copy of $T_{r'+1}$ in $G$ where $x$ plays the
role of the $(s'+1)$th vertex in $T_{r'+1}$. (This is true regardless of whether
$d ^+ _G (x) \geq  ( 1- 1/r  ) n$ or $d ^- _G (x) \geq  ( 1-1/r ) n$.)
This  implies that $G$ contains a
(consistent) copy of $T_{r}$, as desired.
\endproof

The next result is the main tool in the proof of Theorem~\ref{approx}. It will be used to convert a $T_r$-packing in a reduced digraph $R$ into a significantly larger $\{T_r,T_{r+1}\}$-packing.

\begin{lemma}\label{expand}
Let $\eta ,\gamma >0$ and $n,r \geq 2$ be integers such that $0 < 1/n \ll \gamma \ll \eta \ll 1/r$. Let $G$ be a digraph on $n$ vertices with dominant degree sequence $d^* _1 \leq d^* _2 \leq \dots \leq d^* _n$ where
\begin{align}\label{min} d^* _i \geq (r-2)n/r+i +\eta n \ \text{ for all } \ i<n/r.
\end{align}
Further, suppose that the largest $T_r$-packing in $G$ covers precisely $n' \leq (1- \eta )n$ vertices. Then there exists a $\{T_r, T_{r+1} \}$-packing in $G$ that covers at least $n'+ \gamma n$ vertices.
\end{lemma}
\proof
By repeatedly applying Proposition~\ref{fact}, we have that $n' \geq \eta n /2$. 
Define a bijection $I: V(G) \mapsto [n]$ where $I(x)=i$ implies that $d^* _G (x)=d^* _i$. Let $\mathcal M$ denote a $T_r$-packing in $G$ covering precisely $n'$ vertices so that 
$$\sum _{x \in V(G)\setminus V(\mathcal M)} I(x)$$ is maximised.
Set $n'':=n-n'$ and let $G'':=G\setminus V(\mathcal M)$.

\begin{claim}\label{claim}
There are at least $\gamma n$ vertices $x \in V(G'')$ such that $d^* _G (x, V(\mathcal M))\geq (r-1)n'/r +\gamma n$.
\end{claim}
Suppose for a contradiction that the claim is false. Let $V(G'') = \{ x_1 , \dots , x_{n''} \}$ where $I(x_1) < I (x_2) < \dots < I (x_{n''})$. 
For each $ 1\leq i \leq n''$, set $s_i :=I(x_i)$.
So for each $1 \leq i \leq n''$, $d_G ^* (x_i) =d^* _{s_i}$
where $s_1 < s_2 < \dots < s_{n''}$. Let $j$ be the largest integer such that $s_j <n/r$ (so $0\leq j \leq n''$). We call $x_1, \dots , x_j$ \emph{small vertices} and
$x_{j+1}, \dots , x_{n''}$ \emph{big vertices}. 

Notice that (\ref{min}) implies that every big vertex $x \in V(G'')$ satisfies $d^* _G (x) \geq (r-1)n/r +3\eta n/4$.
We say a big vertex $ x \in V(G'')$ is \emph{bad} if $d^* _G (x, V(G'')) \leq (r-1)n''/r +\eta n/2$. 
At most $\gamma n$ big vertices $x \in V(G'')$ are bad.
 Indeed, otherwise there are at least $\gamma n$ vertices $x \in V(G'') $ such that
$$d^* _G (x, V(\mathcal M)) \geq (r-1)n/r +3\eta n/4 -(r-1)n''/r -\eta n/2 \geq (r-1) n'/r +\gamma n,$$ a contradiction to our assumption.

Consider any small vertex $x_i \in V(G'')$. Then $d^* _G (x_i) =d^* _{s_i} \geq (r-2)n/r+s_i +\eta n$ by (\ref{min}) since $s_i < n/r$. Suppose for a contradiction that $d^* _G (x_i, V(G'')) \leq (r-2)n''/r +i+2\gamma n$. Then 
\begin{align}\label{min2}
d^* _G (x_i, V(\mathcal M)) \geq (r-2)n/r +s_i +\eta n -(r-2)n''/r -i-2\gamma n \geq (r-2) n' /r +(s_i -i) +\eta n/2.
\end{align}
Without loss of generality assume that $d^* _G (x_i) =d^+ _G (x_i)$. Then (\ref{min2}) implies that there are at least $(s_i-i)+\eta n/3$ vertices $y \in V(\mathcal M)$ with the property that if $y$ lies in a copy $T'_r$ of
$T_r$ in $\mathcal M$ then $x_i$ sends an edge to every vertex in $V(T'_r) \setminus \{ y\}$. In particular, $(V(T'_r) \cup \{x_i\}) \setminus \{ y\}$ spans a copy of $T_r$ in $G$. Since $I(x_{i'} ) < I(x_i)$ for all $i' <i$, this implies that there is one such vertex $y$ with $s_i=I(x_i) <I(y)$.
Let $\mathcal M '$ denote the $T_r$-packing in $G$ obtained from $\mathcal M$ by replacing the copy $T'_r$ of $T_r$ covering $y$ with the copy of $T_r$ spanning $(V(T'_r) \cup \{x_i\}) \setminus \{ y\}$.
Then $\mathcal M'$ is a $T_r$-packing in $G$ covering precisely $n'$ vertices and with 
$$\sum _{x \in V(G)\setminus V(\mathcal M)} I(x) < \sum _{x \in V(G)\setminus V(\mathcal M')} I(x),$$
a contradiction to the choice of $\mathcal M$. So  $d^* _G (x_i, V(G'')) \geq (r-2)n''/r +i+2\gamma n$.

In summary,  $d^* _G (x_i, V(G'')) \geq (r-2)n''/r +i+2\gamma n$ for each small vertex $x_i \in V(G'')$. Further,  at most $\gamma n$ big vertices $x \in V(G'')$ are bad. By deleting the bad vertices from $G''$ we obtain a subdigraph $G^*$ of $G''$ with dominant degree sequence $d^* _{G^*,1} \leq d^* _{G^*,2} \leq \dots \leq d^* _{G^*, n^*}$ where $n^*:=|G^*|$
so that 
\begin{align*} d^* _{G^*,i} \geq (r-2)n^*/r+i+\gamma n  \ \text{ for all } \ i<n^*/r.
\end{align*}
Hence, Proposition~\ref{fact} implies that $G^*$ contains a copy $T^* _r$ of $T_r$. Then $\mathcal M \cup \{T^* _r\}$ is a $T_r$-packing in $G$ covering more than $n'$ vertices, a contradiction. 
This completes the proof of the claim.

\medskip

Given any $x \in V(G'')$ such that $d^* _G (x, V(\mathcal M))\geq (r-1)n'/r +\gamma n$ there are at least $\gamma n$ copies $T'_r$ of $T_r$ in $\mathcal M$ so that $x$ sends all possible edges to $V(T'_r)$ in
$G$ or  receives all possible edges from $V(T'_r)$ in $G$. In particular, $V(T'_r)\cup \{x\}$ spans a copy of $T_{r+1}$ in  
$G$. Thus, Claim~\ref{claim} implies that  there exists a $\{T_r, T_{r+1} \}$-packing in $G$ that covers at least $n'+ \gamma n$ vertices, as desired.
\endproof

\subsection{Proof of Theorem~\ref{approx}}\label{secapprox}
Define additional constants $\eps, d, \gamma$ and $M' \in \mathbb N$ so that $0<1/n \ll 1/M' \ll \eps \ll d \ll \gamma \ll \eta \ll 1/r, 1/h$. Set $z:= \lceil 1/\gamma \rceil $. Apply Lemma~\ref{dilemma} with parameters $\eps, d$ and $M'$ to $G$ 
to obtain clusters $V_1, \dots , V_k$, an exceptional set $V_0$ and a pure digraph $G'$.
Set $m:=|V_1|=\dots =|V_k|$.
 Let $R$ be the reduced digraph of $G$ with parameters $\eps, d$ and $M'$. Lemma~\ref{inherit} implies that $R$ has 
dominant degree sequence $d^* _{R,1} \leq d^* _{R,2} \leq \dots \leq d^* _{R,k}$ 
such that
\begin{align}\label{ds*} d^* _{R,i} \geq (r-2)k/r+i+\eta k/2  \ \text{ for all } \ i<k/r.
\end{align}

\begin{claim}\label{blowclaim}
 $R':=R(r^z)$ contains a $T_r$-packing covering at least $(1-\eta /2)kr^z=(1-\eta /2)|R'|$ vertices. 
\end{claim}
If $R$ contains a $T_r$-packing covering at least $(1-\eta /2)k$ vertices then Fact~\ref{fact2} implies that Claim~\ref{blowclaim} holds. So suppose that the largest $T_r$-packing in $R$ covers precisely
$d \leq (1- \eta /2)k$ vertices. Then by Lemma~\ref{expand}, $R$ contains a $\{T_r, T_{r+1} \}$-packing that covers at least $d+ \gamma k$ vertices. Thus, by Fact~\ref{fact2}, $R(r)$ contains a $T_r$-packing
covering at least $(d+\gamma k)r$ vertices. (So at least a $\gamma$-proportion of the vertices in $R(r)$ are covered.) Further, Proposition~\ref{fact3} and (\ref{ds*}) imply that 
$R(r)$ has dominant degree sequence $d^* _{R(r),1} \leq \dots \leq d^* _{R(r),kr}$ such that
\begin{align*} d^* _{R(r),i} \geq (r-2)k+i +\eta kr/4 \ \text{ for all } \ i<k.
\end{align*}
If $R(r)$ contains a $T_r$-packing covering at least $(1-\eta /2)kr$ vertices then again Fact~\ref{fact2} implies that the claim holds.
So suppose that the largest $T_r$-packing in $R(r)$ covers precisely
$d' \leq (1- \eta /2)kr$ vertices. Recall that $d' \geq (d+\gamma k)r$. By Lemma~\ref{expand}, $R(r)$ contains a $\{T_r, T_{r+1} \}$-packing that covers at least $d' +\gamma kr \geq (d+2\gamma k)r$ vertices.
Thus, by Fact~\ref{fact2}, $R(r^2)$ contains a $T_r$-packing
covering at least $(d+2\gamma k)r^2$ vertices. (So at least a $2\gamma$-proportion of the vertices in $R(r^2)$ are covered.) Repeating this argument at most $z$ times we see that the claim holds.

\medskip

For each $1 \leq i \leq k$, partition $V_i$ into classes $V^* _i, V_{i,1}, \dots , V_{i, r^z}$ where $m':=|V_{i,j}|= \lfloor m/r^z \rfloor \geq m/(2r^z)$ for all $1 \leq j \leq r^z$.
Since $mk \geq (1-\eps)n$ by Lemma~\ref{dilemma},
\begin{align}\label{m'}
m'|R'| = \big \lfloor {m}/{r^z} \big \rfloor kr^z \geq mk-kr^z \geq (1-2\eps)n.
\end{align}

Lemma~\ref{slice} implies that if $(V_{i_1}, V_{i_2})_{G'}$ is $\eps$-regular with density at least $d$ then $(V_{i_1, j_1}, V_{i_2,j_2})_{G'}$ is $2\eps r^z$-regular with density at least $d-\eps \geq d/2$
(for all $1\leq j_1,j_2 \leq r^z$). In particular, we can label the vertex set of $R'$ so that $V(R')=\{V_{i,j} : 1 \leq i \leq k , \ 1 \leq j \leq r^z \}$ where 
$V_{i_1, j_1} V_{i_2,j_2} \in E(R')$ implies that $(V_{i_1, j_1}, V_{i_2,j_2})_{G'}$ is $2\eps r^z$-regular with density at least $d/2$. 

By Claim~\ref{blowclaim}, $R'$ has a $T_r$-packing $\mathcal M$ that contains at least $(1-\eta /2)|R'|$ vertices. Consider any copy $T'_r$ of $T_r$ in $\mathcal M$ and let
$V(T'_r)=\{ V_{i_1, j_1}, V_{i_2,j_2}, \dots , V_{i_r, j_r} \}$. Set $V'$ to be the union of $ V_{i_1, j_1}, V_{i_2,j_2}, \dots , V_{i_r, j_r}$. Note that $0<1/m' \ll 2 \eps r^z \ll d/2 \ll \gamma \ll 1/r, 1/h$.
Thus, Lemma~\ref{red} implies that $G'[V']$ contains a $T^h _r$-packing covering all but at most $\sqrt{2 \eps r ^z} m' r \leq \gamma m' r$ vertices. By considering each copy of $T_r$ in $\mathcal M$
we conclude that
$G'\subseteq G$ contains a $T^h _r$-packing covering at least
$$ (1-\gamma )m'r \times (1- \eta /2)|R'|/r \stackrel{(\ref{m'})}{\geq} (1-\gamma)(1-\eta /2) (1-2\eps )n \geq (1-\eta )n $$
vertices, as desired.
\endproof

\section{Overview of our  absorbing approach}\label{overabs}
The `absorbing method' was first used in~\cite{rrs2} and has subsequently
been  applied to numerous embedding problems in extremal graph theory. In this section we give an overview of the absorbing method that we will apply to prove Theorems~\ref{absthm} and~\ref{absorbG}.

Suppose that in some (di)graph (or hypergraph) $G$ we wish to find a `small' $H$-absorbing set $M$ in $G$ that absorbs any `very small' set of vertices in $G$. To prove that such an $H$-absorbing set $M$
exists it is sufficient to show that $G$ has many `connecting structures' of a certain type. That is,  it suffices to show that  there
is a constant $t \in \mathbb N$ such that, given any $x,y \in V(G)$, there are `many' $t$-sets $X \subseteq V(G)$ so that both $G[X \cup \{x\}]$ and $G[X \cup \{y\}]$ contain perfect $H$-packings.
The following lemma of Lo and Markstr\"om~\cite{lo} makes this observation precise.

\begin{lemma}[Lo and Markstr\"om~\cite{lo}]\label{lo}
Let $h,t \in \mathbb N$ and let $\gamma >0$. Suppose that $H$ is a hypergraph on $h$ vertices. Then there exists an $n_0 \in \mathbb N$ such that the following holds. Suppose that $G$ is a hypergraph
on $n \geq n_0$ vertices so that, for any $x,y \in V(G)$, there are at least $\gamma n^{th-1}$ $(th-1)$-sets $X \subseteq V(G)$ such that both $G[X \cup \{x\}]$ and $G[X \cup \{y\}]$ contain perfect $H$-packings.
Then $V(G)$ contains a set $M$ so that
\begin{itemize}
\item $|M|\leq (\gamma/2)^h n/4$;
\item $M$ is an $H$-absorbing set for any $W \subseteq V(G) \setminus M$ such that $|W| \in h \mathbb N$ and  $|W|\leq (\gamma /2)^{2h} hn/32 $.
\end{itemize}
Moreover, the analogous result holds for digraphs $H$ and $G$.
\end{lemma} 
(Lo and Markstr\"om~\cite{lo} only proved Lemma~\ref{lo} in the case of hypergraphs. However, the proof of the digraph case is identical.)

In view of Lemma~\ref{lo} it is natural to seek sufficient conditions that guarantee our desired connecting structures $X$ for each $x,y \in V(G)$. For this, we introduce the definition of an \emph{$H$-path} below. In what follows we will be working in the (di)graph setting, although everything naturally extends to hypergraphs.

Suppose that $H$ is a (di)graph on $h$ vertices. An \emph{$H$-path} $P$ is a (di)graph with the following properties:
\begin{itemize}
\item[(i)] $|P|=th+1$ for some  $t \in \mathbb N$;
\item[(ii)]  $V(P)=X_1 \cup \dots \cup X_t \cup \{y_1, \dots , y_{t+1}\}$ where $|X_i|=h-1$ for all $1 \leq i \leq t$;
\item[(iii)]  $P[ X_i \cup \{y_i\}] =P[ X_i \cup \{y_{i+1}\}]=H$ for all $1\leq i \leq t$.
\end{itemize}
We call $y_1$ and $y_{t+1}$ the \emph{endpoints of $P$} and $t$ the \emph{length of $P$}. 
An example of a $T_3$-path of length $3$ is given in Figure~2. 
 A \emph{truncated $H$-path} $Q$ is a (di)graph with the following properties:
\begin{itemize}
\item[(i)] $|Q|=th-1$ for some  $t \in \mathbb N$;
\item[(ii)]  $V(Q)=X_1 \cup \dots \cup X_t \cup \{y_2, \dots , y_{t}\}$ where $|X_i|=h-1$ for all $1 \leq i \leq t$;
\item[(iii)]  $Q[ X_i \cup \{y_i\}] =Q[ X_i \cup \{y_{i+1}\}]=H$ for all $2\leq i \leq t-1$ and $Q[X_1\cup \{y_2\}]=Q[X_t\cup \{y_t\}]=H$.
\end{itemize}
We call $X_1$ and $X_t$ the \emph{endsets of $Q$} and $t$ the \emph{length of $Q$}. Notice that given any $H$-path $P$ of length $t$ and endpoints $x$ and $y$, $Q:=P\setminus\{x,y\}$ is a truncated $H$-path
of length $t$. In this case we say that $Q$ is \emph{the truncated $H$-path of $P$}.

The following simple observation follows immediately from the definition of an $H$-path.

\begin{figure}\label{picture}
\begin{center}\footnotesize
\includegraphics[width=0.5\columnwidth]{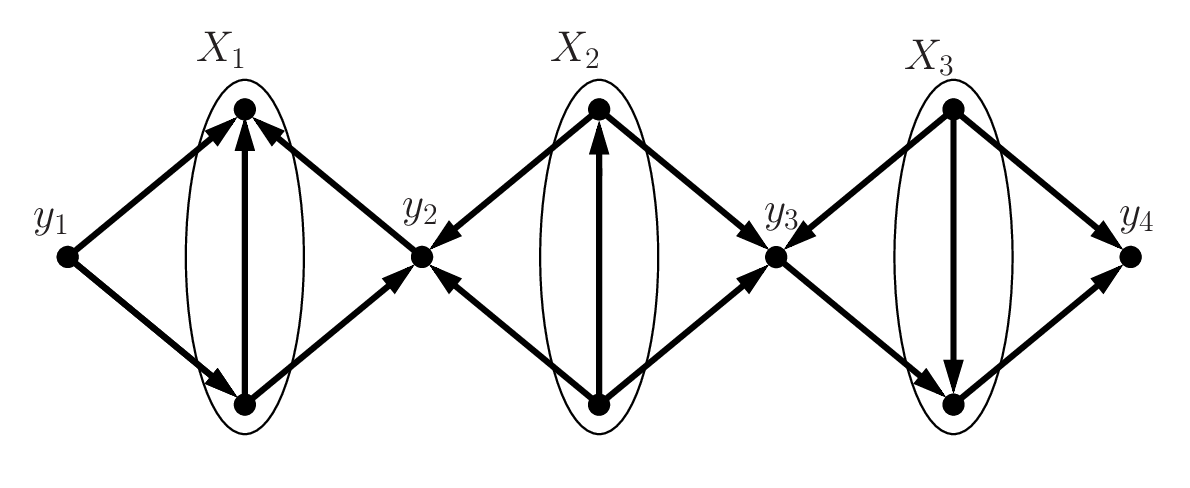}
\caption{An example of a $T_3$-path of length $3$}
\end{center}
\end{figure}

\begin{fact}\label{Hfact} Let $H$ be a (di)graph.
 Let $P_1$ be an $H$-path of length $t_1$ with endpoints $x,y$ and let $P_2$ be an $H$-path of length $t_2$ with endpoints $y,z$. Suppose that $P_1 \setminus \{y\}$ and $P_2 \setminus \{y\}$ are vertex-disjoint. Then $P_1 \cup P_2$ is an $H$-path of length $t_1+t_2$ with endpoints $x,z$.
\end{fact}

Suppose that $P$ is an $H$-path in a (di)graph $G$ with endpoints $x$ and $y$. Set $X:= V(P) \setminus \{x,y\}$. Then by condition (iii) in the definition of an $H$-path,  both $G[X \cup \{x\}]$ and $G[X \cup \{y\}]$ contain perfect $H$-packings. Thus, to show that the hypothesis of Lemma~\ref{lo} holds it suffices to prove that, given any $x,y \in V(G)$, there are at least $\gamma n^{th-1}$ $(th-1)$-sets $X \subseteq V(H)$ so that $X \cup \{x,y\}$ spans an $H$-path in $G$ of length $t$ with endpoints $x$ and $y$. We will use this approach in the proofs of Theorems~\ref{absthm} and~\ref{absorbG}.

Let $0< \gamma \ll \beta \ll 1/t,1/h$ where $t,h \in \mathbb N$. Suppose that $H$ is a (di)graph on $h$ vertices and $G$ is a sufficiently large (di)graph on $n$ vertices. Define an auxiliary graph $\mathcal G$ as follows.
$\mathcal G$ has vertex set $V(G)$ and given distinct $x,y \in V(G)$, $xy \in E(\mathcal G)$ precisely when there are at least $\beta n^{h-1}$ $H$-paths in $G$ of length one with endpoints $x$ and $y$.
Suppose that for each $x,y \in V(G)$, there are at least $\beta n^{t-1}$ paths of length $t$ in $\mathcal G$ with endpoints $x$ and $y$. Then Fact~\ref{Hfact} together with a simple calculation implies that there are at least $\gamma n^{th-1}$
$H$-paths in $G$ of length $t$ with endpoints $x$ and $y$. Thus, ultimately to find the  connecting structures required to apply Lemma~\ref{lo} it suffices to prove that the auxiliary graph $\mathcal G$ is `well-connected'. Although $\mathcal G$ will not be formally required,  it
 will serve as a useful point of reference when describing how we construct our desired $T_r$-paths in the proof of Theorem~\ref{absthm}.



\section{The absorbing results for Theorems~\ref{conj2} and~\ref{mainthm2}}\label{sec8}
\subsection{Proof of Theorem~\ref{absthm}}\label{secab1}
In this section we prove Theorem~\ref{absthm}. For some fixed $t'\in \mathbb N$, we will show that, for any $x,y \in V(G)$, there are `many' $T_r$-paths of length $t'$ in $G$ with endpoints $x$ and $y$.
Then by applying Lemma~\ref{lo} this implies that $G$ contains our desired absorbing set $M$.

Let $G$ be as in Theorem~\ref{absthm}, $H=T_r$ and recall the definition of the auxiliary graph $\mathcal G$ given in Section~\ref{overabs} (where the constant $\beta$ in the definition of $\mathcal G$ is chosen
such that $\beta \ll \eta , 1/r$). 
Roughly speaking, the next result states that given any vertex $x \in V(\mathcal G)$, there are `many' vertices $y$ such that $xy \in E(\mathcal G)$ and $y$  has `large' dominant degree (relative to $x$). This `expansion' property will then be exploited in Lemma~\ref{abs2} to conclude that $\mathcal G$ is
`well-connected', and thus $G$ contains the desired $T_r$-paths.

\begin{lemma}\label{abs1}
Let $\gamma, \alpha , \eta >0$ and $n,r \geq 2$ be integers such that $0 < 1/n \ll \gamma \ll \alpha \ll \eta \ll 1/r$. Let $G$ be a digraph on $n$ vertices with dominant degree sequence $d^* _1 \leq d^* _2 \leq \dots \leq d^* _n$ where
\begin{align}\label{gen} d^* _i \geq (r-2)n/r+i +\eta n \ \text{ for all } \ i<n/r.
\end{align}
Suppose that $x \in V(G)$ such that 
\begin{align}\label{gen1} d^*_G (x) \geq (r-2)n/r+j +\eta n \ \text{ for some } \ 1 \leq j<n/r-\alpha n.
\end{align}
Then there exist at least $\gamma n ^{r}$ pairs $(y,X)$ where $y \in V(G)\setminus \{x\}$, $X \subseteq V(G)\setminus \{x,y\}$ and the following conditions hold:
\begin{itemize}
\item[(i)] $d^*_G (y) \geq (r-2)n/r+j +(\eta+\gamma) n  $;
\item[(ii)] $X \cup \{x\}$ and $X \cup \{y\}$ both spans copies of $T_r$ in $G$. In particular, $X\cup\{x,y\}$ spans a copy of a $T_r$-path of length one in $G$ with endpoints $x$ and $y$.
\end{itemize}
\end{lemma}
\proof
Without loss of generality assume that $d^+ _G (x) \geq (r-2)n/r +j +\eta n$.
Suppose first that $r=2$. Then (\ref{gen}) and (\ref{gen1}) imply that there are at least $\alpha n$ vertices $z \in N^+ _G (x)$ such that $d^* _G (z) \geq j+(\eta +\alpha /2)n$. Fix such a vertex $z$ and
without loss of generality suppose that $d^* _G (z)=d^+ _G (z)$. Then there are at least $\alpha n$ vertices $y \in N^+ _G (z) \setminus \{x\}$ such that $d^* _G (y) \geq j+(\eta +\gamma)n$. Set $X: = \{z\}$.
Then $(y,X)$ satisfies (i) and (ii). Further, there are at least
$$\alpha n \times \alpha n  > \gamma n^2$$ 
choices for $(y,X)$, as desired.

Next suppose that $r \geq 3$.
 Let $r' \in \mathbb N$. Suppose that $T'$ is a copy of $T_{r'}$ in $G$.
Let $V(T')=\{x_1, \dots , x_{r'}\}$ where $x_i$ plays the role of the $i$th vertex
of $T_{r'}$. We say that $T'$ is \emph{consistent} if there exists $0\leq s' \leq r'$
such that
\begin{itemize}
\item $d_G ^+ (x_i) \geq  ( 1- {1}/{r} + 3 \eta /4 ) n$ for all $i \leq s'$;
\item $d_G ^- (x_i) \geq  ( 1- {1}/r + 3 \eta /4) n$ for all $i >s'$.
\end{itemize}
We call $s'$ a \emph{turning point of $T'$}.

\begin{claim}\label{c19}
There are at least $\eta ^{r-1} n^{r-2}$ sets $X' \subseteq N^+ _G (x)$ such that $X'$ spans a consistent copy of $T_{r-2}$ in $G$.
\end{claim}
By (\ref{gen}) and (\ref{gen1}) there are at least 
$(r-2){n}/{r} +j +\eta n -{n}/{r} \geq \eta n$ vertices $y_1 \in N^+ _G (x)$ such that $d^* _G (y_1) \geq (1-1/r+3 \eta /4)n$.
In particular, $\{y_1\}$ is a consistent copy of $T_1$. So if $r=3$ the claim follows. Thus, assume $r \geq 4$.

Suppose for some $1 \leq r' \leq r-3$, $T'$ is a consistent copy of $T_{r'}$ in $G$ such that $V(T') \subseteq N^+ _G (x)$.
Let $V(T')=\{y_1, \dots , y_{r'}\}$ where $y_i$ plays the role of the $i$th vertex
of $T_{r'}$ and let $s'$ denote a turning point of $T'$. 
Set $N':=N^+ _G (x) \cap \bigcap _{i\leq s'} N^+ _G(y_i) \cap \bigcap _{i>s'} N^- _G (y_i).$
Since $T'$ is consistent with turning point $s'$ and $r'\leq r-3$, (\ref{gen1}) implies that
$$|N'|\geq (r-2)n/r+j +\eta n - r'n/r \geq n/r +\eta n .$$
In particular, (\ref{gen}) implies that there are at least $\eta n$ vertices $w \in  N'$ such that $d_G ^* (w) \geq ( 1- {1}/{r} + 3 \eta /4 ) n$. Notice that 
$V(T') \cup \{w\}\subseteq N^+ _G (x)$ spans a consistent copy of $T_{r'+1}$ in $G$ where $w$ plays the
role of the $(s'+1)$th vertex in $T_{r'+1}$. (This is true regardless of whether $d_G ^+ (w) \geq ( 1- {1}/{r} + 3 \eta /4 ) n$ or $d_G ^-(w) \geq ( 1- {1}/{r} + 3 \eta /4 ) n$.)

This argument shows that there are at least 
$$(\eta n)^{r-2} \times \frac{1}{(r-2)!} \geq \eta ^{r-1} n^{r-2}$$ sets $X' \subseteq N^+ _G (x)$ such that $X'$ spans a consistent copy of $T_{r-2}$ in $G$, thus proving the claim.

\medskip

Fix a set $X'$ as in Claim~\ref{c19}. Let $X'=\{ x_1, \dots , x_{r-2} \}$ where $x_i$ plays the role of the $i$th vertex of a consistent copy $T'_{r-2}$ of $T_{r-2}$ in $G$. 
Let $s$ denote a turning point of $T'_{r-2}$. Set $N:= N^+ _G (x) \cap \bigcap _{i\leq s} N^+ _G(x_i) \cap \bigcap _{i>s} N^- _G (x_i).$
Then $$|N| \stackrel{(\ref{gen1})}{\geq} (r-2)n/r +j +\eta n -(r-2)n/r =j+\eta n.$$
In particular, (\ref{gen}) and (\ref{gen1}) imply that there are at least $\alpha n$ vertices $z \in N$ such that $d^* _G (z) \geq (r-2)n/r+j +(\eta +\alpha /2 )n$.

Fix such a vertex $z$ and let $X:=X' \cup \{z\}$. Note that $X \cup \{x\}$ spans a copy of $T_r$ in $G$. Without loss of generality suppose that $d^+ _G (z)  \geq (r-2)n/r+j +(\eta +\alpha /2 )n$.
Set $N^* := N^+ _G (z) \cap \bigcap _{i\leq s} N^+ _G(x_i) \cap \bigcap _{i>s} N^- _G (x_i).$
Then $|N^*| \geq (r-2)n/r +j +(\eta +\alpha /2)n-(r-2)n/r =j +(\eta +\alpha /2)n$.
Thus, by (\ref{gen}) there are at least $\alpha n $ vertices $y \in N^* \setminus \{x\}$ such that $d_G ^* (y) \geq  (r-2)n/r+j +(\eta+\gamma) n $. Further, $X \cup \{y\}$ spans a copy of $T_r$ in $G$ (this follows  by the choice of $z$ and $y$).  Fix such a vertex $y$. So (i) and (ii) hold.

Note that there were at least $\eta ^{r-1} n^{r-2}$ choices for $X'$, $\alpha n$ choices for $z$ and $\alpha n$ choices for $y$. So in total there are at least
$$\eta ^{r-1} n^{r-2} \times \alpha n \times \alpha n \times \frac{1}{r-1} \geq \gamma n^r$$ choices for $(y,X)$, as desired.
\endproof

The next result guarantees our desired $T_r$-paths in $G$ (implicitly implying that our auxiliary graph $\mathcal G$ is `well-connected').

\begin{lemma}\label{abs2}
Let $\xi, \gamma, \alpha , \eta >0$ and $n,r \geq 2$ be integers such that $0 < 1/n \ll \xi  \ll \gamma \ll \alpha \ll \eta \ll 1/r$. Let $G$ be a digraph on $n$ vertices with dominant degree sequence $d^* _1 \leq d^* _2 \leq \dots \leq d^* _n$ where
\begin{align}\label{gen2} d^* _i \geq (r-2)n/r+i +\eta n \ \text{ for all } \ i<n/r.
\end{align}
Set $t:=(2\lfloor 1/\gamma r -\alpha/\gamma \rfloor +1)r-1$ and $t':=2\lfloor 1/\gamma r -\alpha/\gamma \rfloor +1$. Given any distinct $x,y \in V(G)$ there are at least $\xi n^t$ $t$-sets $A \subseteq V(G) \setminus \{x,y\}$ so that $A\cup \{x,y\}$ spans a $T_r$-path of length $t'$ in $G$ with endpoints $x$ and $y$.
In particular, for each such set $A$,
$G[A \cup \{x\}]$ and $G[A \cup \{y\}]$ contain perfect $T_r$-packings.
\end{lemma}
\proof
Define an additional constant $\beta $ so that $\xi \ll \beta \ll \gamma$.
Set $t'':= \lfloor 1/\gamma r -\alpha/\gamma \rfloor$.
\begin{claim}\label{claimcon}
Let $1 \leq i \leq t''$. There are at least $ \beta ^i n^{ir}$ pairs $(w,X)$ where $w \in V(G) \setminus \{x,y \}$, $X \subseteq V(G) \setminus \{x,y,w \}$ and $|X|=ir-1$ with the following properties:
\begin{itemize}
\item $d^*_G (w) \geq (r-2)n/r +(\eta+i \gamma) n  $;
\item $X \cup \{x,w\}$ spans a $T_r$-path of length $i$ in $G$ with endpoints $x$ and $w$. In particular,
$G[X \cup \{x\}]$ and $G[X \cup \{w\}]$ both contain perfect $T_r$-packings.
\end{itemize}
\end{claim}
We prove the claim inductively.
By (\ref{gen2}), $d^* _G (x) \geq (r-2)n/r+1+\eta n$. Thus, by Lemma~\ref{abs1}, Claim~\ref{claimcon} holds for $i=1$. Suppose that for some $1\leq i <t''$, Claim~\ref{claimcon} holds.
Let $(w,X)$ be as in Claim~\ref{claimcon} with $|X|=ir-1$. Set $j:=i\gamma n < n/r-\alpha n$. Note that 
$$d^* _G (w) \geq (r-2)n/r + (\eta+i \gamma) n = (r-2)n/r + j +\eta n.$$ 
By Lemma~\ref{abs1} there exist at least $\gamma n^r$ pairs $(w',X')$ where $w' \in V(G) \setminus \{w\}$, $X' \subseteq V(G)\setminus \{w,w'\}$ and the following conditions hold:
\begin{itemize}
\item $d^*_G (w') \geq (r-2)n/r+j +(\eta+\gamma) n =(r-2)n/r+(\eta +(i+1)\gamma)n $;
\item $X' \cup \{w,w'\}$ spans a $T_r$-path of length one in $G$ with endpoints $w$ and $w'$. That is,
$X' \cup \{w\}$ and $X' \cup \{w'\}$ both spans copies of $T_r$ in $G$.
\end{itemize}
Fix such a pair $(w',X')$ so that $X' \cup \{w'\}$ is disjoint from $X \cup \{x,y \}$; There are at least
$$\gamma n^r -r(ir+1) \binom{n}{r-1} \geq \gamma n^r /2$$
choices for $(w',X')$.

Set $X'':=X \cup \{w\} \cup X'$. Then by Fact~\ref{Hfact}, $(w', X'')$ is such that $w' \in V(G) \setminus \{x,y\}$, $X'' \subseteq V(G) \setminus \{x,y,w'\}$ and $|X''|=(i+1)r-1$ with the following properties:
\begin{itemize}
\item $d^*_G (w') \geq (r-2)n/r +(\eta+(i+1) \gamma) n  $;
\item $X'' \cup \{x,w'\}$ spans a $T_r$-path of length $i+1$ in $G$ with endpoints $x$ and $w'$. In particular,
$G[X'' \cup \{x\}]$ and $G[X'' \cup \{w'\}]$ both contain perfect $T_r$-packings.
\end{itemize}
There are at least $\beta ^i n^{ir}$ choices for $(w,X)$ and at least $\gamma n^r /2$ choices for $(w',X')$. Thus, in total there are at least 
$$\beta ^{i}  n^{ir} \times \frac{\gamma n^r }{2} \times \frac{1}{((i+1)r)!} \geq \beta ^{i+1} n^{(i+1)r}$$
such pairs $(w',X'')$. Thus, Claim~\ref{claimcon} holds with respect to $i+1$, as desired. 

\medskip

Fix a pair $(w,X)$ as in Claim~\ref{claimcon} with $i:=t''$. The proof of the next claim is analogous to that of Claim~\ref{claimcon} (so we omit it).
\begin{claim}\label{claimcon2}
Let $1 \leq i \leq t''$. There are at least $ \beta ^i n^{ir}$ pairs $(z,Y)$ where $z \in V(G)  $, $Y \subseteq V(G)\setminus \{z\} $ and $|Y|=ir-1$ with the following properties:
\begin{itemize}
\item $Y\cup \{z\}$ is disjoint from $X\cup \{w,x,y\}$;
\item $d^*_G (z) \geq (r-2)n/r +(\eta+i \gamma) n  $;
\item  $Y \cup \{y,z\}$ spans a $T_r$-path of length $i$ in $G$ with endpoints $y$ and $z$. In particular, $G[Y \cup \{y\}]$ and $G[Y \cup \{z\}]$ both contain perfect $T_r$-packings.
\end{itemize}
\end{claim}
Fix a pair $(z,Y)$ as in Claim~\ref{claimcon2} with $i:=t''$. Note that 
\begin{align}\label{star}
d^* _G (w), d^* _G (z) \geq (r-2)n/r +(\eta +t'' \gamma )n \geq (r-1)n/r +(\eta -\alpha -\gamma)n \geq (r-1)n/r +3\eta n/4.
\end{align}
Without loss of generality assume that $d^+ _G (w)=d_G ^* (w)$ and $d^+ _G (z)=d_G ^* (z)$ (the other cases  follow analogously). Recall the definition of a consistent copy of $T_{r'}$ in $G$ given in Lemma~\ref{abs1}. By applying (\ref{star}) and arguing in a similar way to the proof of Claim~\ref{c19} we obtain the following claim.

\begin{claim}\label{claimcon3} Suppose that $r\geq 3$. Then 
there are at least $\eta ^{r-1} n^{r-2}$ sets $Z' \subseteq N^+ _G (w)\cap N^+ _G (z)$ such that $Z'$ spans a consistent copy of $T_{r-2}$ in $G$.
\end{claim}

Suppose that $r \geq 3$.
Fix a set $Z'$ as in Claim~\ref{claimcon3} so that $Z'$ is disjoint from $X\cup Y \cup \{x,y\}$; There are at least 
$$\eta ^{r-1} n^{r-2} - 2t''r\binom{n}{r-3} \geq \eta ^{r-1} n^{r-2} /2$$ choices for $Z'$. Let $Z'=\{ z_1, \dots , z_{r-2} \}$ where $z_i$ plays the role of the $i$th vertex of a consistent copy $T'_{r-2}$ of $T_{r-2}$ in $G$. 
Let $s$ denote a turning point of $T'_{r-2}$. Set $N:= N^+ _G (w) \cap N^+_G (z) \cap \bigcap _{i\leq s} N^+ _G(z_i) \cap \bigcap _{i>s} N^- _G (z_i).$
Then $$|N| \stackrel{(\ref{star})}{\geq} (r-2)n/r  +3\eta n/2 -(r-2)n/r =3\eta n /2.$$
Fix some $z' \in N\setminus (X\cup Y \cup \{x,y\})$. Notice that there are at least $\eta n$ choices for $z'$. Set $Z:=Z' \cup \{z'\}$. The choice of $z'$ implies that $Z \cup \{w\}$ and $Z\cup \{z\}$ span copies of
$T_r$ in $G$.

If $r=2$, then notice that there are at least $\eta n$ vertices $z' \in N^+ _G (w) \cap N^+_G (z)$ that are not in $X \cup Y \cup \{x,y\}$. Fix such a $z'$ and let $Z:=\{z'\}$.

For every $r \geq 2$, set $A:= X\cup Y \cup Z \cup \{w,z\}$. So $|A|=(2t''+1)r-1=t$. By construction of $A$, $A \cup \{x,y\}$ spans a $T_r$-path of length $t'$ in $G$ with endpoints $x$ and $y$.
If $r \geq 3$, then recall that there are at least $\beta ^{t''} n^{t''r}$ choices for $(w,X)$, at least $\beta ^{t''} n^{t''r}$ choices for $(z,Y)$, at least $\eta ^{r-1} n^{r-2}/2$ choices for $Z'$ and at least $\eta n$ choices for $z'$.
Overall, this implies that there are at least
$$ \beta ^{t''} n^{t''r} \times \beta ^{t''} n^{t''r}\times \eta ^{r-1} n^{r-2}/2 \times \eta n \times \frac{1}{t!} \geq \xi n^{t}$$ choices for $A$.
If $r=2$, then a similar calculation shows there are at least $\xi n^{t}$ choices for $A$, as required.
\endproof

We now deduce Theorem~\ref{absthm} from Lemmas~\ref{lo} and~\ref{abs2}.

\noindent
{\bf Proof of Theorem~\ref{absthm}.}
Define additional constants $\xi , \gamma , \alpha$ such that $\nu \ll \xi \ll \gamma \ll \alpha \ll \eta$. Set $t:=(2\lfloor 1/\gamma r -\alpha/\gamma \rfloor +1)r-1$ and $t':=2\lfloor 1/\gamma r -\alpha/\gamma \rfloor +1$ and define $\gamma ' >0$ so that $\nu= (\gamma '/2)^r /4$. Let $G$ be as in the statement of the theorem.  Lemma~\ref{abs2} implies that, for any $x,y \in V(G)$, there are at least
$\xi n^t$ $t$-sets $X \subseteq V(G) \setminus \{x,y\}$ such that both $G[X \cup \{x\}]$ and $G[X \cup \{y\}]$ contain perfect $T_r$-packings.
Since $\nu \ll \xi \ll 1/r$, and by definition of $\gamma '$, we have that $\xi n^t \geq \gamma ' n^t = \gamma 'n^{t'r-1}$.
Apply Lemma~\ref{lo} with $T_r, r,t',\gamma '$ playing the roles of $H,h,t,\gamma$. 
Then $V(G)$ contains a set $M$ so that
\begin{itemize}
\item $|M|\leq (\gamma'/2)^r n/4=\nu n$;
\item $M$ is a $T_r$-absorbing set for any $W \subseteq V(G) \setminus M$ such that $|W| \in r \mathbb N$ and  $|W|\leq \nu ^3 n \leq (\gamma '/2)^{2r} rn/32 $.
\end{itemize}
\endproof

\subsection{Proof of Theorem~\ref{absorbG}}\label{secabs2}
In this section we prove Theorem~\ref{absorbG}. We will show that there is some fixed $s'\in \mathbb N$ such that, for any $x,y \in V(G)$, there are `many' $K^h _r$-paths $P$ of length $s'$ in $G$ with endpoints $x$ and $y$. Consider such a $K^h _r$-path $P$ and set $X:=V(P) \setminus \{x,y\}$.
Since $K^ h_ r$ contains a perfect $H$-packing (for any $H$ such that $|H|=h$ and $\chi (H)=r$), the definition of a $K^h _r$-path implies that   $G[X\cup \{x\}]$ and $G[X \cup \{y\}]$ contain perfect $H$-packings.
Then by applying Lemma~\ref{lo} this implies that $G$ contains our desired $H$-absorbing set $M$.

The following result guarantees our desired $K^h _r$-paths in the case when $h=1$. 
It is an immediate consequence of Lemma~\ref{abs2}.
\begin{lemma}\label{abs3}
Let $\xi, \gamma, \alpha , \eta >0$ and $n,r \geq 2$ be integers such that $0 < 1/n \ll \xi  \ll \gamma \ll \alpha \ll \eta \ll 1/r$. Let $G$ be a graph on $n$ vertices with  degree sequence $d _1 \leq d _2 \leq \dots \leq d _n$ where
\begin{align*} d _i \geq (r-2)n/r+i +\eta n \ \text{ for all } \ i<n/r.
\end{align*}
Set $t:=(2\lfloor 1/\gamma r -\alpha/\gamma \rfloor +1)r-1$ and $t':=2\lfloor 1/\gamma r -\alpha/\gamma \rfloor +1$. Given any distinct $x,y \in V(G)$ there are at least $\xi n^t$ $t$-sets $A \subseteq V(G) \setminus \{x,y\}$ so that $A\cup \{x,y\}$ spans a $K_r$-path of length $t'$ in $G$ with endpoints $x$ and $y$.
In particular, for each such set $A$,
$G[A \cup \{x\}]$ and $G[A \cup \{y\}]$ contain perfect $K_r$-packings.
\end{lemma}
We will apply Lemma~\ref{abs3} to a reduced graph $R$ of $G$ in order to conclude that  $R$ contains a  $K_r$-path between any two clusters. We then show that this $K_r$-path in $R$ corresponds to  `many' $K^h_r$-paths in $G$. For this, we first show that a special type of  `blow-up'
of a $K_r$-path is itself a $K^h _r$-path.

\medskip

Let $h, t,r \in \mathbb N$ such that $r+1, t \geq 3$. Consider a $K_r$-path $P$ of length $t$. Let $X_1,\dots, X_t, y_1,\dots , y_{t+1}$ be as in the definition of  $P$. In particular, $P$ has endpoints
$y_1$ and $y_{t+1}$ and $|X_i|=r-1$ for all $1 \leq i \leq t$. We define $P^*(h)$ to be the graph obtained from $P$ as follows:
\begin{itemize}
\item Replace each vertex set $X_i$ with a set of $h(r-1)$ vertices $X_i (h)$ so that $X_i (h)$ induces a copy of $K^h _{r-1}$  in $P^* (h)$ (for $1 \leq i \leq t$);
\item For $2\leq i \leq t-1$, we replace $y_i$ with a set $Y_i$ of $h$ vertices; $y_t$ is replaced by a set $Y_t$ of $2h-1$ vertices. Set $Y_1:=\{ y_1\}$, $Y_{t+1} := \{y_{t+1}\}$;
\item In $P^* (h)$ there are all possible edges between $Y_i \cup Y_{i+1}$ and $X_i (h)$ (for all $1\leq i \leq t$).
\end{itemize}
Given a truncated $K_r$-path $Q$ of length $t$, we define $Q^*(h)$ analogously. In particular, if $Q$ is the truncated $K_r$-path of $P$ then $Q^* (h) =P^* (h) \setminus \{y_1,y_{t+1}\}$.

Note that $P^* (h)$ is a blow-up of $P$ where the vertices in $V(P) \setminus \{y_1, y_t,y_{t+1}\}$ are replaced by $h$ vertices, $y_t$ is replaced by $2h-1$ vertices and $y_1$ and $y_{t+1}$ are left `untouched'.
In particular, this immediately implies the following fact.

\begin{fact}\label{blowfact2}
Let $h, t,r \in \mathbb N$ such that $r+1, t \geq 3$ and let $G$ and $H$ be graphs. Suppose that $P$ is a $K_r$-path of length $t$  and $Q$ is a truncated $K_r$-path of length $t$. If $P\subseteq G$ then $P^* (h) \subseteq G(2h-1)$. Further, if 
$Q \subseteq H$ then $Q^* (h)  \subseteq H(2h-1)$.
\end{fact}

We now show that, crucially, $P^* (h)$ is a $K^h _r$-path of length $t$.
\begin{lemma}\label{star2}
Let $h, t,r \in \mathbb N$ such that $r+1, t \geq 3$. 
\begin{itemize}
\item[(a)]
Suppose that $P$ is a $K_r$-path of length $t$ with endpoints $y_1$ and $y_{t+1}$. Then $P^* (h)$ is a $K_r ^h$-path of length $t$ with endpoints $y_1$ and $y_{t+1}$. 
\item[(b)]
Suppose that $Q$ is a  truncated $K_r$-path  of length $t$ with endsets $X_1$ and $X_{t}$. Then $Q^* (h)$ is a truncated $K_r ^h$-path of length $t$. Further, let $Q'$ denote the graph obtained from $Q^* (h)$
by adding a vertex $x$ that is adjacent to every vertex in $X_1 (h)$ and a vertex $y$ that is adjacent to every vertex in $X_t (h)$. Then $Q'$ is a $K^h _r$-path of length $t$ with endpoints $x$ and $y$.
\end{itemize}
\end{lemma}
\proof We first prove (a).
 Let $X_1,\dots, X_t, y_1,\dots , y_{t+1}$ be as in the definition of  $P$ and $X_1 (h),\dots,$ $X_t(h)$, $Y_1,\dots , Y_{t+1}$ be as in the definition of $P^* (h)$.
Firstly, note that
$$|P^* (h)|= \left ( \sum _{1 \leq i \leq t } |X_i (h)| \right )+ 2 + \left ( \sum _{2 \leq i \leq t-1 } |Y_i| \right ) +|Y_{t}| = h(r-1)t +2 + h(t-2) +(2h-1) =hrt+1.$$
So $P^* (h)$ satisfies condition (i) of the definition of a $K^h _r$-path of length $t$.

\smallskip 

For each $2 \leq i \leq t-1$, define $Y' _{i-1}$, $y'_{i}$ so that $|Y' _i|=h-1$ and $Y_i = Y' _{i-1} \cup \{ y'_{i}\}$. Define $Y'_{t-1}$, $Y'_{t}$, $y'_t$ so that $|Y' _{t-1}|=|Y'_t|=h-1$ and 
$Y_t = Y' _{t-1} \cup Y' _{t} \cup \{ y'_{t}\}$. Set $X'_i := X_i (h) \cup Y'_i$ for all $1 \leq i \leq t$ and $y'_1:=y_1$, $y'_{t+1}:=y_{t+1}$. Thus, $V(P^* (h))=X'_1 \cup \dots \cup X'_t \cup \{y'_1, y'_2, \dots  ,y'_{t+1}\}$ where
$|X'_i|=h(r-1)+h-1= hr-1$ for all $1 \leq i \leq t$. So $P^* (h)$ satisfies condition (ii) of the definition of a $K^h _r$-path of length $t$.

\smallskip

Let $1 \leq i \leq t-1$. Recall that $X_i (h)$ induces a copy of $K^h _{r-1}$ in $P^* (h)$. Since $Y'_i \subseteq Y_{i+1}$ and there are all possible edges between $Y_{i+1}$ and $X_i (h)$ in $P^* (h)$, we have that
$X'_i$ induces a copy of the complete $r$-partite graph $K^*$ with $r-1$ vertex classes of size $h$ and one vertex class (namely $Y'_i$) of size $h-1$. A similar argument implies that $X'_t$ induces a copy of $K^*$ in
$P^* (h)$. 
Recall that $y'_i \in Y_i$ for all $1 \leq i \leq t+1$. So by definition of $P^* (h)$, $y'_i$ and $y'_{i+1}$ send an edge to every vertex in $X_i (h)$ in $P^* (h)$ (for all $1 \leq i \leq t$). Altogether, this
implies that $X'_i \cup \{y'_i\}$ and $X'_i \cup \{y'_{i+1}\}$ induce copies of $K^h _r$ in $P^* (h)$. Thus, condition (iii) of the definition of a $K^h _r$-path of length $t$ is satisfied. So (a) holds.
It is easy to see that (b) follows from (a).
\endproof

We now apply Lemmas~\ref{abs3} and~\ref{star2} to conclude that a graph $G$ as in Theorem~\ref{absorbG} contains many $K^h _r$-paths of a given length between any distinct $x,y \in V(G)$.

\begin{lemma}\label{abs9}
Let $\gamma ', \gamma, \alpha , \eta >0$ and $n,r ,h \in \mathbb N$ such that $0 < 1/n \ll \gamma '  \ll \gamma \ll \alpha \ll \eta \ll 1/r , 1/h$ where $r \geq 2$. Let $G$ be a graph on $n$ vertices with  degree sequence $d _1 \leq d_2 \leq \dots \leq d _n$ where
\begin{align}\label{gen9} d _i \geq (r-2)n/r+i +\eta n \ \text{ for all } \ i<n/r.
\end{align}
Set  $s':=2\lfloor 1/\gamma r -\alpha/\gamma \rfloor +3$ and $s:=hrs'-1$. Given any distinct $x,y \in V(G)$ there are at least $\gamma ' n^s$ $s$-sets $A \subseteq V(G) \setminus \{x,y\}$ so that $A\cup \{x,y\}$ spans a $K^h _r$-path of length $s'$ in $G$ with endpoints $x$ and $y$.
In particular, for each such set $A$,
$G[A \cup \{x\}]$ and $G[A \cup \{y\}]$ contain perfect $K^h _r$-packings.
\end{lemma}
\proof
 Let $\eps, d, \gamma , \alpha , \eta >0$ and $M' \in \mathbb N$ such that $$0<1/M' \ll \eps \ll d \ll \gamma \ll \alpha \ll \eta \ll 1/r,1/h.$$
Set 
$t':=2\lfloor 1/\gamma r -\alpha/\gamma \rfloor +1$. Let $M$ denote the integer obtained by applying Lemma~\ref{dilemma} with parameters
$\eps$ and $M'$. Let $ \gamma '>0$ and $n \in \mathbb N$ be sufficiently large such that
$$0<1/n  \ll \gamma ' \ll 1/M', 1/M.$$

Suppose that $G$ is a graph on $n$ vertices as in the statement of the lemma. Apply Lemma~\ref{dilemma} with parameters $\eps, d$ and $M'$ to $G$ to obtain clusters 
$V_1, \dots , V_k$, an exceptional set $V_0$ and a pure graph $G'$.
Set $m:=|V_1|=\dots =|V_k|$. So $M' \leq k \leq M$ and $(1-\eps)n/M \leq m \leq n/M'$.
 Let $R$ denote the reduced graph of $G$ with parameters $\eps, d$ and $M'$. Lemma~\ref{inherit2} implies that $R$ has 
 degree sequence $d_{R,1} \leq d _{R,2} \leq \dots \leq d _{R,k}$ 
such that
\begin{align}\label{dsa} d _{R,i} \geq (r-2)k/r+i+\eta k/2  \ \text{ for all } \ i<k/r.
\end{align}

Consider distinct $x,y \in V(G)$. 
 We may assume that $x,y \in V_0$. (Otherwise, we delete the at most two clusters containing $x$ or $y$, and move
their vertices into $V_0$. The same properties of $R$ still hold, just with slightly perturbed parameters $k$ and $\eta$.)
Let $N_R (x)$ denote the set of clusters $V_i \in V(R)$ such that $d_G (x,V_i)\geq \alpha m$. Set $d_R (x) := |N_R (x)|$. Define $N_R (y)$ and $d_R (y)$ analogously.
Note that 
$$ m d_R (x) + k\alpha m +|V_0| \geq d_G (x) \stackrel{(\ref{gen9})}{\geq} (r-2)n/r +\eta n.$$
So as $km\leq n$ and $|V_0| \leq \eps n$,
\begin{align}\label{XR}
d_R (x) \geq (r-2)k/r+\eta k/2.
\end{align}
Similarly, $d_R (y) \geq (r-2)k/r+\eta k/2.$

Using (\ref{dsa}) and (\ref{XR}), a simple greedy argument implies that there is a copy $K_x$ of $K_{r-1}$ in $R$ so that:
\begin{itemize}
\item $V(K_x) \subseteq N_R (x)$;
\item $V(K_x)=\{U_1, \dots , U_{r-1}\}$ where $d_R (U_i)\geq (r-1)k/r+\eta k/3$ for all $1 \leq i \leq r-2$ and $d_R (U_{r-1}) \geq (r-2)k/r+\eta k/2$.
\end{itemize}

Note that the latter condition implies that there is a cluster $V_x \in V(R)$ such that $V(K_x) \cup \{V_x\}$ spans a copy of $K_r$ in $R$.
Similarly, there is a copy $K_y$ of $K_{r-1}$ in $R$ and a cluster $V_y \in V(R)$ so that:
\begin{itemize}
\item $V(K_y) \subseteq N_R (y) $;
\item $V(K_y)=\{W_1, \dots , W_{r-1}\}$ where $d_R (W_i)\geq (r-1)k/r+\eta k/3$ for all $1 \leq i \leq r-2$ and $d_R (W_{r-1}) \geq (r-2)k/r+\eta k/2$;
\item $V(K_y) \cup \{V_y\}$ is vertex-disjoint from $V(K_x) \cup \{V_x\}$ and spans a copy of $K_r$ in $R$.
\end{itemize}

Remove the clusters in $K_x$ and $K_y$ from $R$ and denote the resulting graph by $R'$. Set $k':=|R'|$. (\ref{dsa}) implies that  $R'$ has 
 degree sequence $d_{R',1} \leq d _{R',2} \leq \dots \leq d _{R',k'}$ 
such that
\begin{align}\label{dsb} d _{R',i} \geq (r-2)k'/r+i+\eta k'/3  \ \text{ for all } \ i<k'/r.
\end{align}

Apply Lemma~\ref{abs3} with $R', \eta /3$ playing the roles of $G, \eta$. So $R'$ contains a $K_r$-path $P$ of length $t'$  with endpoints $V_x$ and $V_y$. (Here, we do not need to use the fact that
there are `many' such $K_r$-paths.) By the choice of $V_x$ and $ V_y$, $\mathcal V:=V(P)\cup V(K_x) \cup V(K_y)$ spans a truncated $K_r$-path $Q$ of length $t'+2=s'$ in $R$  with endsets $X_1:=V(K_x)$ and $X_{s'}:=
V(K_y)$.

For each $V \in \mathcal V$ let $V'$ denote a subset of $V$ of size $\alpha m$ so that:
\begin{itemize}
\item $V'\subseteq N_G (x)$ if $V \in V(K_x)$;
\item $V'\subseteq N_G (y)$ if $V \in V(K_y)$.
\end{itemize}
(So for $V \in V(P)$ the choice of $V'$ is arbitrary.)
Notice that we can choose such subsets $V'$ since $V(K_x) \subseteq N_R (x)$ and $V(K_y) \subseteq N_R (y)$.
Suppose that $V,W \in \mathcal V$ and $VW \in E(R)$. Then since $(V,W)_{G'}$ is $\eps$-regular with density at least $d$, Lemma~\ref{slice} implies that $(V',W')_{G'}$ is $\eps ^{1/2}$-regular with density at least $d/2$.

Set $R'':=R[\mathcal V]$.
Since $Q\subseteq R''$, Fact~\ref{blowfact2} implies that $Q_1:=Q^* (h)$ is a subgraph of $R''(2h-1)$. Moreover, Lemma~\ref{star2}(b) implies that $Q_1$ is a truncated $K^h_r$-path of length $s'$.
Let $G''$ denote the subgraph of $G$ whose vertex set consists of all vertices $v$ such that $v \in V'$ for some $V \in \mathcal V$ and so that $G''[V',W']=G'[V',W']$ for all $V,W \in \mathcal V$.
Apply Lemma~\ref{simo} with $G'',R'', Q_1 ,2h-1, \eps ^{1/2}, d/2, \alpha m$  playing the roles of $G,R,H, t, \eps , d, m$. Note that $|Q_1|=s$ and set $\Delta :=\Delta (Q_1) $ and $\eps _0 := (d/2-\eps ^{1/2})^\Delta /(2+\Delta)$. As $\eps \leq \eps _0$ and $2h-2 \leq \eps _0 \alpha m$, Lemma~\ref{simo} implies that there are at least
$$(\eps _0 \alpha m)^s \geq (d^{\Delta+1} \alpha (1-\eps )n/M)^s \geq s! \gamma ' n^s$$ 
copies of $Q_1 $ in $G'' \subseteq G$. 

Consider any copy $Q'_1$ of $Q_1$ guaranteed by Lemma~\ref{simo}. Recall that $Q$ has endsets $X_1 =V(K_x) \subseteq N_R (x)$ and $X_{s'} =V(K_y) \subseteq N_R (y)$. Lemma~\ref{simo} guarantees that the vertices in the set $X_1(h)\subseteq V(Q'_1)$ corresponding to the $r-1$ elements of $X_1$ are embedded into sets $V'\subseteq N_G (x)$. Similarly, the vertices in the set $X_{s'} (h)\subseteq V(Q'_1)$ corresponding to the $r-1$ elements of $X_{s'}$ are embedded into sets $V'\subseteq N_G (y)$. Thus, Lemma~\ref{star2}(b) implies that $V(Q'_1) \cup \{x,y\}$ spans a $K^h _r$-path in $G$ of length
$s'$  with endpoints $x$ and $y$.

Since there are at least $s! \gamma ' n^s$ such copies $Q'_1$ of $Q_1$ in $G''$, there are at least $\gamma ' n^s$ $s$-sets $A \subseteq V(G)$ such that $A=V(Q'_1)$ for some such $Q'_1$. 
So indeed, there are at least $\gamma ' n^s$ $s$-sets $A \subseteq V(G) \setminus \{x,y\}$ so that $A\cup \{x,y\}$ spans a $K^h _r$-path of length $s'$ in $G$ with endpoints $x$ and $y$.
\endproof

We now deduce Theorem~\ref{absorbG} from Lemmas~\ref{lo} and~\ref{abs9}.

\noindent
{\bf Proof of Theorem~\ref{absorbG}.}
Define additional constants $ \gamma , \alpha$ such that $\nu  \ll \gamma \ll \alpha \ll \eta$. Set
$s':=2\lfloor 1/\gamma r -\alpha/\gamma \rfloor +3$ and $s:=hrs'-1$ and define $\gamma ' >0$ so that $\nu= (\gamma '/2)^h /4$.
In particular, $\gamma ' \ll \gamma$.

 Let $G$ be as in the statement of the theorem.  Lemma~\ref{abs9} implies that, for any $x,y \in V(G)$, there are at least
$\gamma' n^s$ $s$-sets $X \subseteq V(G) \setminus \{x,y\}$ such that both $G[X \cup \{x\}]$ and $G[X \cup \{y\}]$ contain perfect $K^h _r$-packings and thus, perfect $H$-packings.
Apply Lemma~\ref{lo} with $rs',\gamma '$ playing the roles of $t,\gamma$. 
Then $V(G)$ contains a set $M$ so that
\begin{itemize}
\item $|M|\leq (\gamma'/2)^h n/4=\nu n$;
\item $M$ is an $H$-absorbing set for any $W \subseteq V(G) \setminus M$ such that $|W| \in h \mathbb N$ and  $|W|\leq \nu ^3 n \leq (\gamma '/2)^{2h} rn/32 $.
\end{itemize}
\endproof

\section*{Acknowledgements}
The author would like to thank Peter Keevash for a helpful conversation, in particular, for suggesting the auxiliary graph $\mathcal G$ in Section~\ref{overabs} as a tool for presenting the absorbing method.
The author is also grateful to Daniela K\"uhn and Deryk Osthus  for helpful comments, and to the referees for the careful reviews.

{\footnotesize \obeylines \parindent=0pt
Andrew Treglown
School of Mathematics
University of Birmingham
Edgbaston
Birmingham
B15 2TT
UK
\tt{a.c.treglown@bham.ac.uk}}
\end{document}